\documentclass[a4paper, 12pt]{amsart}

\usepackage{amsmath,amssymb,enumitem,verbatim,stmaryrd,xcolor,microtype}
\usepackage[T1]{fontenc}
\usepackage[utf8]{inputenc}
\usepackage[english]{babel}
\usepackage[headheight=14pt,top=3.3cm,bottom=3.3cm,left=3.2cm,right=3.2cm]{geometry}
\usepackage[linktoc=page,colorlinks,linkcolor={red!80!black},citecolor={red!80!black},urlcolor={blue!80!black}]{hyperref}\newcommand{\arxiv}[1]{\href{http://arxiv.org/pdf/#1}{arXiv:#1}}
\usepackage[all]{xy}
\usepackage[backgroundcolor=yellow,linecolor=yellow,textsize=footnotesize]{todonotes} \makeatletter \providecommand \@dotsep{5} \def\listtodoname{List of Todos} \def\listoftodos{\@starttoc{tdo}\listtodoname} \makeatother 
\usepackage{tikz,calc,tikz-cd}
\usetikzlibrary{matrix,arrows,decorations.pathreplacing}
\pgfrealjobname{overviewF1}
\usepackage{pgfplots}
\pgfplotsset{width=7cm,compat=1.13}

\usepackage{etoolbox}
\makeatletter
\patchcmd{\@startsection}{\@afterindenttrue}{\@afterindentfalse}{}{}\makeatother    
\patchcmd{\section}{\scshape}{\bfseries}{}{}\makeatletter\renewcommand{\@secnumfont}{} 
\patchcmd{\@settitle}{\uppercasenonmath\@title}{\large}{}{}
\patchcmd{\@setauthors}{\MakeUppercase}{}{}{}
\addto{\captionsenglish}{} 
\addto{\captionsenglish}{} 
\addto{\captionsenglish}{} 
\makeatother

\usepackage{fancyhdr}

\theoremstyle{plain}
\newtheorem{thm}{Theorem}[section]

\newtheorem{lemma}[thm]{Lemma}

\newtheorem*{question*}{Question}

\theoremstyle{definition}

\newtheorem{rem}[thm]{Remark}
\newtheorem{ex}[thm]{Example}
\newtheorem*{RH}{Riemann hypothesis}

\usepackage{mathptmx}                                        
\DeclareRobustCommand{\gobblefour}[4]{}    
\DeclareSymbolFont{sfoperators}{OT1}{bch}{m}{n} \DeclareSymbolFontAlphabet{\mathsf}{sfoperators} \makeatletter\def\operator@font{\mathgroup\symsfoperators}\makeatother 
\DeclareSymbolFont{cmletters}{OML}{cmm}{m}{it}              
\DeclareSymbolFont{cmsymbols}{OMS}{cmsy}{m}{n}
\DeclareSymbolFont{cmlargesymbols}{OMX}{cmex}{m}{n}
\DeclareMathSymbol{\myjmath}{\mathord}{cmletters}{"7C}     \let\jmath\myjmath 
\DeclareMathSymbol{\myamalg}{\mathbin}{cmsymbols}{"71}     \let\amalg\myamalg
\DeclareMathSymbol{\mycoprod}{\mathop}{cmlargesymbols}{"60}\let\coprod\mycoprod
\DeclareMathSymbol{\myalpha}{\mathord}{cmletters}{"0B}     \let\alpha\myalpha 
\DeclareMathSymbol{\mybeta}{\mathord}{cmletters}{"0C}      \let\beta\mybeta
\DeclareMathSymbol{\mygamma}{\mathord}{cmletters}{"0D}     \let\gamma\mygamma
\DeclareMathSymbol{\mydelta}{\mathord}{cmletters}{"0E}     \let\delta\mydelta
\DeclareMathSymbol{\myepsilon}{\mathord}{cmletters}{"0F}   \let\epsilon\myepsilon
\DeclareMathSymbol{\myzeta}{\mathord}{cmletters}{"10}      \let\zeta\myzeta
\DeclareMathSymbol{\myeta}{\mathord}{cmletters}{"11}       \let\eta\myeta
\DeclareMathSymbol{\mytheta}{\mathord}{cmletters}{"12}     \let\theta\mytheta
\DeclareMathSymbol{\myiota}{\mathord}{cmletters}{"13}      \let\iota\myiota
\DeclareMathSymbol{\mykappa}{\mathord}{cmletters}{"14}     \let\kappa\mykappa
\DeclareMathSymbol{\mylambda}{\mathord}{cmletters}{"15}    \let\lambda\mylambda
\DeclareMathSymbol{\mymu}{\mathord}{cmletters}{"16}        \let\mu\mymu
\DeclareMathSymbol{\mynu}{\mathord}{cmletters}{"17}        \let\nu\mynu
\DeclareMathSymbol{\myxi}{\mathord}{cmletters}{"18}        \let\xi\myxi
\DeclareMathSymbol{\mypi}{\mathord}{cmletters}{"19}        \let\pi\mypi
\DeclareMathSymbol{\myrho}{\mathord}{cmletters}{"1A}       \let\rho\myrho
\DeclareMathSymbol{\mysigma}{\mathord}{cmletters}{"1B}     \let\sigma\mysigma
\DeclareMathSymbol{\mytau}{\mathord}{cmletters}{"1C}       \let\tau\mytau
\DeclareMathSymbol{\myupsilon}{\mathord}{cmletters}{"1D}   \let\upsilon\myupsilon
\DeclareMathSymbol{\myphi}{\mathord}{cmletters}{"1E}       \let\phi\myphi
\DeclareMathSymbol{\mychi}{\mathord}{cmletters}{"1F}       \let\chi\mychi
\DeclareMathSymbol{\mypsi}{\mathord}{cmletters}{"20}       \let\psi\mypsi
\DeclareMathSymbol{\myomega}{\mathord}{cmletters}{"21}     \let\omega\myomega
\DeclareMathSymbol{\myvarepsilon}{\mathord}{cmletters}{"22}\let\varepsilon\myvarepsilon
\DeclareMathSymbol{\myvartheta}{\mathord}{cmletters}{"23}  \let\vartheta\myvartheta
\DeclareMathSymbol{\myvarpi}{\mathord}{cmletters}{"24}     \let\varpi\myvarpi
\DeclareMathSymbol{\myvarrho}{\mathord}{cmletters}{"25}    \let\varrho\myvarrho
\DeclareMathSymbol{\myvarsigma}{\mathord}{cmletters}{"26}  \let\varsigma\myvarsigma
\DeclareMathSymbol{\myvarphi}{\mathord}{cmletters}{"27}    \let\varphi\myvarphi

\DeclareMathOperator{\Spec}{Spec}

\DeclareMathOperator{\Sch}{Sch}

\DeclareMathOperator{\Hom}{Hom}

\DeclareMathOperator{\Sets}{Sets}

\DeclareMathOperator{\Mat}{{Mat}}
\DeclareMathOperator{\GL}{{GL}}
\DeclareMathOperator{\SL}{{SL}}
\DeclareMathOperator{\PGL}{{PGL}}
\DeclareMathOperator{\Gr}{{Gr}}
\DeclareMathOperator{\Stab}{{Stab}}
\DeclareMathOperator{\Norm}{{Norm}}
\DeclareMathOperator{\Cent}{{Cent}}
\DeclareMathOperator{\TSch}{{Sch_{\mathcal{T}}}}

\newcommand\A{{\mathbb A}}
\newcommand\B{{\mathbb B}}
\newcommand\C{{\mathbb C}}
\newcommand\F{{\mathbb F}}
\newcommand\G{{\mathbb G}}

\newcommand\N{{\mathbb N}}
\renewcommand\P{{\mathbb P}}
\newcommand\Q{{\mathbb Q}}
\newcommand\R{{\mathbb R}}

\newcommand\T{{\mathbb T}}
\newcommand\Z{{\mathbb Z}}

\newcommand\cC{{\mathcal C}}

\newcommand\cG{{\mathcal G}}

\newcommand\cM{{\mathcal M}}

\newcommand\cO{{\mathcal O}}
\newcommand\cP{{\mathcal P}}

\newcommand\cR{{\mathcal R}}

\newcommand\cT{{\mathcal T}}
\newcommand\cU{{\mathcal U}}

\newcommand\cW{{\mathcal W}}
\newcommand\cX{{\mathcal X}}

\newcommand\fe{{\mathfrak e}}
\newcommand\fm{{\mathfrak m}}
\newcommand\fp{{\mathfrak p}}

\renewcommand\Re{{\textup{Re}\;}}
\newcommand\Fun{{\F_1}}
\newcommand\Funsq{{\F_{1^2}}}
\newcommand\Funn{{\F_{1^n}}}

\newcommand\res{\textup{res}}
\newcommand\id{\textup{id}}

\newcommand\trop{\textup{trop}}

\newcommand\rk{\textup{rk}}
\newcommand\barx{\overline{x}}

\newcommand\sT{{\scriptscriptstyle\cT\hspace{-3pt}}}

\newcommand\overZ{\overline{\Spec\Z}}   

\renewcommand\={\equiv}
\renewcommand\geq{\geqslant}
\renewcommand\leq{\leqslant}
\newcommand{\gen}[1]{\langle #1 \rangle}
\newcommand{\norm}[1]{|#1|}
\newcommand{\bpquot}[2]{#1\!\sslash\!#2}
\renewcommand{\binom}[2]{\Bigl(\begin{matrix}#1\\#2\end{matrix}\Bigr)}
\newcommand{\qbinom}[2]{\Bigl[\,\begin{matrix}#1\\#2\end{matrix}\,\Bigr]_q}

\newcommand{\bpgenquot}[2]{#1\!\sslash\!\gen{#2}}
\newdir{ >}{{}*!/-5pt/@^{(}}\newcommand{\arincl}[1]{\ar@{ >->}@<-0,0ex>#1} 
\newcommand\tvdots{\raisebox{3pt}{$\scalebox{.75}{\vdots}$}}

\newenvironment{psmallmatrix}{\left(\begin{smallmatrix}}{\end{smallmatrix}\right)}

\title{$\Fun$ for everyone}

\author{Oliver Lorscheid}
\email{oliver@impa.br}
\address{Instituto Nacional de Matem\'atica Pura e Aplicada, Rio de Janeiro, Brazil}

\begin{document}

\begin{abstract}
 This text serves as an introduction to $\Fun$-geometry for the general mathematician. We explain the initial motivations for $\Fun$-geometry in detail, provide an overview of the different approaches to $\Fun$ and describe the main achievements of the field.
\end{abstract}

\maketitle



\section*{Prologue}
\label{prologue}

$\Fun$-geometry is a recent area of mathematics that emerged from certain heuristics in combinatorics, number theory and homotopy theory that could not be explained in the frame work of Grothendieck's scheme theory. These ideas involve a theory of algebraic geometry over a hypothetical field $\Fun$ with one element, which includes a theory of algebraic groups and their homogeneous spaces, an arithmetic curve $\overZ$ over $\Fun$ that compactifies the arithmetic line $\Spec\Z$ and $K$-theory that stays in close relation to homotopy theory of topological spaces such as spheres.

Starting in the mid 2000's, there was an outburst of different approaches towards such theories, which bent and extended scheme theory. Meanwhile some of the initial problems for $\Fun$-geometry have been settled and new applications have been found, for instance in cyclic homology and tropical geometry.

The first thought that crosses one's mind in this context is probably the question:
\medskip
\begin{center}\it 
 What is the ``field with one element''? 
\end{center}
\medskip
Obviously, this oxymoron cannot be taken literally as it would imply a mathematical contradiction. It is resolved as follows. First of all we remark that we do not need to define $\Fun$ itself---what is needed for the aims of $\Fun$-geometry is a suitable category of schemes over $\Fun$.

However, many approaches contain an explicit definition of $\Fun$, and in most cases, the field with one element is not a field and has two elements. Namely, the common answer of many theories is that $\Fun$ is the multiplicative monoid $\{0,1\}$, lacking any additive structure.

Continuing the word plays, the protagonist $\Fun$ is also written as $\F_{\textup{un}}$, coming from the French word ``\emph{un}'' for $1$, and it enters the title of the paper ``\emph{Fun with $\Fun$}'' of Connes, Consani and Marcolli (\cite{Connes-Consani-Marcolli09}). In other sources, $\Fun$-geometry disguises itself as ``\emph{non-additive geometry}'' or ``\emph{absolute geometry}''.


\subsection*{Content overview}
The first chapter \ref{section: conception} of this text is dedicated to a description of the main ideas that led to $\Fun$-geometry. We explain the two leading problems in detail: Tits' dream of explaining certain combinatorial geometries in terms of an algebraic geometry over $\Fun$ in sections \ref{subsection: incidence geometry} and the ambitious programme to prove the Riemann hypothesis in section \ref{subsection: Riemann hypothesis}.

The second chapter \ref{section: birth hour} serves as an overview of the manifold approaches towards $\Fun$-geometry. We will describe two prominent theories of $\Fun$-schemes in detail: Deitmar's theory of monoid schemes in section \ref{subsection: monoid schemes} and the author's theory of blue schemes in section \ref{subsection: blueprints}.

The third chapter \ref{section: growing up} summarizes the impact of $\Fun$-geometry on mathematics today. We spend a few words on developments around the Riemann hypothesis in section \ref{subsection: steps towards the Riemann hypothesis}, describe in detail the realization of Tits' dream via blue schemes in section \ref{subsection: algebraic groups over F1} and outline promising recent developments in tropical geometry that involve $\Fun$-schemes in section \ref{subsection: applications to tropical geometry}.


\section{Conception: the heuristics leading to \texorpdfstring{$\Fun$}{F1}}\label{section: conception}

The first mentioning of a ``field of characteristic one'' in the literature can be found in the 1957 paper \cite{Tits57} by Jacques Tits. The postulation of such a field origins in his observation that certain projective geometries over finite fields $\F_q$ with $q$ elements have a meaningful analogue for $q=1$. Tits remarks that these latter geometries should have an explanation in terms of projective geometry over a field with one element. 

In the modern literature on the topic, the analogue for $q=1$ is often called the ``limit $q\to 1$'', a notion that finds its origin in the connections to quantum groups where $q$ occurs indeed as a complex parameter.

For instance, the group of invertible matrices with coefficient in $\F_q$ converges towards the symmetric group $S_n$ on $n$ elements,
\[
 \GL(n,\F_q) \quad \xrightarrow[q\to1]{} \quad \GL(n,\Fun) \ = \ S_n,
\]
compatible with the respective actions on Grassmannians $\Gr(k,n)(\F_q)$ and the family of $k$-subsets of $\{1,\dotsc,n\}$. This is explained in detail in section \ref{subsection: incidence geometry}.

According to a private conversation with Cartier, Tits' idea did not find much resonance at the time---one has to bear in mind that this was at a moment in which the community struggled with a generalization of algebraic geometry from fields to other types of rings; Grothendieck's clarifying invention of schemes was still several years ahead. In so far a geometry over the hypothetical object $\Fun$ was too far away from conceivable mathematics at the time. 

As a consequence, it took more than three decades until the field with one element gained popularity, this time thanks to one of the most famous riddles in number theory, the Riemann hypothesis. Alexander Smirnov gave talks about how $\Fun$-geometry could be involved in a proof of the Riemann hypothesis in the late 1980s. This idea finds its first mentioning in the literature in Manin's influential lecture notes \cite{Manin95}, which are based on his talks at Harvard, Yale, Columbia and MSRI in 1991/92. 

In a nutshell, this ansatz postulates a completed arithmetical curve $\overZ$ over $\Fun$, which would allow one to mimic Weil's proof for function fields. A particular ingredient of this line of thought is that $\Z$ is an algebra over $\Fun$ and that there is a base extension functor
\[
 -\otimes_\Fun\Z: \ \bigl\{ \, \Fun\text{-algebras / schemes}\,\bigr\} \quad \longrightarrow \quad \bigl\{ \, \Z\text{-algebras / schemes}\,\bigr\}.
\]
This is explained in more detail in section \ref{subsection: Riemann hypothesis}.

Around the same time, Smirnov explained another possible application of $\Fun$-geometry in \cite{Smirnov92}: conjectural Hurwitz inequalities for a hypothetical map 
\[
 \overZ \quad \longrightarrow \quad \P^1_\Fun
\]
from the completed arithmetical curve to the projective line over $\Fun$ would imply the abc-conjecture.

Soon after, Kapranov and Smirnov aim in the unpublished note \cite{Kapranov-Smirnov95} to calculate cohomological invariants of arithmetic curves in terms of cohomology over $\Funn$. The unfinished text contains an outburst of different ideas: linear and homological algebra over $\Funn$, distinguished morphisms as cofibrations, fibrations and equivalences (which might be seen as a first hint of the connections of $\Fun$-geometry to homotopy theory), Arakelov theory {modulo $n$} and connections to class field theory and reciprocity laws, which can be seen in analogy to knots and links in $3$-space.

In \cite{Soule04}, Soul\'e explains a connection to the stable homotopy groups of spheres, an idea that he attributes to Manin. Elaborating the formula $\GL(n,\Fun)=S_n$, there should be isomorphisms 
\[
 K_*(\Fun) \quad = \quad \pi_*({B\GL(\infty,\Fun)}^+) \quad  = \quad \pi_*({BS_\infty}^+) \quad \simeq  \quad \pi_*^s(\mathbb{S}) 
\]
where the first equality is the definition of $K$-theory via Quillen's plus construction, naively applied to the elusive field $\Fun$. The second equality is derived from the hypothetical formula 
\[
 \GL(\infty,\Fun) \quad = \quad \bigcup_{n\geq 1} \GL(n,\Fun) \quad = \quad \bigcup_{n\geq 1} S_n \quad = \quad S_\infty.
\]
The last isomorphism is the Barratt-Priddy-Quillen theorem.


\subsection{Incidence geometry and \texorpdfstring{$\Fun$}{F1}} \label{subsection: incidence geometry}
In his seminal paper \cite{Tits57} from 1957, Tits investigates analogues of homogeneous spaces for Lie groups over finite fields. In the following, we explain his ideas in the example of the general linear group $\GL(n)$ of invertible $n\times n$-matrices acting by base change on a Grassmannian $\Gr(k,n)$ of $k$-dimensional subspaces of an $n$-dimensional vector space.\footnote{To be precise, Tits considers in \cite{Tits57} only semi-simple algebraic groups and he considers $\PGL(n)$ in place of $\GL(n)$. However, we can illustrate Tits' idea in the case of either group and we will allow ourselves this inaccuracy for the sake of a simplified account.}

We remark that all these thoughts should find a much more conceptual explanation within the theory of buildings that was introduced a few years later Tits; the interested reader will find more information on the developments of buildings in \cite{Rousseau09}. However, we will refrain from such a reformulation in order to stay historically accurate and to avoid burdening the reader with an introduction to buildings.

\subsubsection{Tits' notion of a geometry}
The example of $\GL(n)$ acting on $\Gr(k,n)$ makes sense over the real or complex numbers as well as over a finite field $\F_q$ with $q$ elements. In the latter case, we are concerned with the group $G=\GL(n,\F_q)$ and the $\F_q$-rational points $\Gr(k,\F_q^n)=\Gr(k,n)(\F_q)$ of the Grassmannian. 

Since the action of $G$ on $\Gr(k,\F_q^n)$ is transitive, $\Gr(k,\F_q^n)$ stays in bijection to the left cosets of the stabilizer of a $k$-subspace $V$ of $\F_q^n$ in $G$. If we choose $V$ to be spanned by the first $k$ standard basis vectors, then the stabilizer $P_k=\Stab_G(V)$ consists of all matrices in $G$ of the form
\[
 \begin{psmallmatrix}
  \begin{smallmatrix}
     \ast  & \dotsb &  \ast  \\[-3pt]
   \tvdots &        &\tvdots \\[-2pt]
     \ast  & \dotsb &  \ast 
  \end{smallmatrix} &
  \begin{smallmatrix}
     \ast  & \dotsb &  \ast  \\[-3pt]
   \tvdots &        &\tvdots \\[-2pt]
     \ast  & \dotsb &  \ast 
  \end{smallmatrix} \\
  \begin{smallmatrix}
      0    & \dotsb &   0    \\[-3pt]
   \tvdots &        &\tvdots \\[-2pt]
      0    & \dotsb &   0   
  \end{smallmatrix}  &
  \begin{smallmatrix}
     \ast  & \dotsb &  \ast  \\[-3pt]
   \tvdots &        &\tvdots \\[-2pt]
     \ast  & \dotsb &  \ast 
  \end{smallmatrix}
 \end{psmallmatrix}
\]
where the upper left block contains an invertible $k\times k$-matrix and the lower right block contains an invertible $(n-k)\times (n-k)$-matrix. Thus we obtain an identification of $\Gr(k,\F_q^n)$ with $G/P_k$.

The containment relation $V'\subset V$ of different subspaces of $\F_q^n$ defines an incidence relation $\iota$ between the elements $V'\in\Gr(k',\F_q^n)$ and $V\in\Gr(k,\F_q^n)$ for different $k'$ and $k$. 

Tits dubs a collection of various homogeneous spaces for a fixed group $G$ together with an incidence relation a \emph{geometry}. He investigates various properties that are satisfied by geometries coming from matrix groups over $\F_q$, like the one described above or its analogues for symplectic groups or orthogonal groups.

The name ``geometry'' can be motivated in the above example. The points of the different Grassmannians $\Gr(1,\F_q^n)$, $\Gr(2,\F_q^n)$, $\dotsc$, $\Gr(n-1,\F_q^n)$ correspond to the points, lines, $\dotsc$, $(n-2)$-dimensional subspaces of the projective space $\P^{n-1}(\F_q)=\Gr(1,\F_q^n)$ and the action of $G$ on the different Grassmannians corresponds to the permutation of linear subspaces by the action of $G$ on $\P^{n-1}(\F_q)$.

\subsubsection{The limit geometry}
Note that every finite field $\F_q$ with $q$ elements produces such a geometry. In other words, the geometry depends on the ``parameter'' $q$. The crucial observation that led Tits to postulate the existence of a field $\Fun$ with $1$ element is that there is a meaningful limit of a geometry when $q$ goes to $1$. More precisely, for every geometry coming from a matrix group over $\F_q$, there is a geometry satisfying the same aforementioned properties and which looks like the limit $q\to 1$.

We explain this limit in our example of $G=\GL(n,\F_q)$ and the Grassmannians $\Gr(k,\F_q^n)$ for various $n$. The group for the limit geometry is the symmetric group $S_n$ on $n$ elements. The homogeneous spaces are the families $\Sigma(k,n)$ of all $k$-subsets of $\{1,\dotsc,n\}$. The incidence relation is defined in terms of the containment $X'\subset X$ for different subsets $X'$ and $X$ of $\{1,\dotsc,n\}$. Note that the stabilizer of the subset $\{1,\dotsc,k\}$ is $S_k\times S_{n-k}$, thus we have $\Sigma(k,n)=S_n/(S_k\times S_{n-k})$.

As mentioned before, the geometry that consists of the homogeneous spaces $\Sigma(k,n)$ of $S_n$ satisfies analogous properties to the geometry of Grassmannians $\Gr(k,\F_q^n)$. 

A first link between these two geometries is laid in terms of the Weyl groups of $\GL(n)$ and the stabilizers $P_k$. To explain, the Weyl group of a matrix group $G$ is defined as the quotient $W=\Norm_G(T)/T$ of the normalizer of the diagonal torus $T$ in $G$ by $T$ itself.\footnote{This is, again, slightly inaccurate. In general, one can consider the Weyl group for any torus of a matrix group. However, if the torus is not specified, it is assumed that the torus is of maximal rank. For $G=\GL(n)$, the diagonal torus is of maximal rank, but this is not true for all matrix groups.} In our case, the normalizer of $T$ consists of all monomial matrices, i.e.\ matrices that have precisely one non-zero entry in each row and each column. The elements of the quotient $W=\Norm_G(T)/T$ can thus be represented by permutation matrices and we conclude that the Weyl group $W$ of $G$ is isomorphic to $S_n$. Similarly, the Weyl group of $P_k$ is $S_k\times S_{n-k}$.

The idea that the geometry of the $\Sigma(k,n)$ should be thought as the limit $q\to 1$ of the geometry of the $\Gr(k,\F_q^n)$ is suggested by the behaviour of the invariants counting points, lines, et cetera.

To start with, these counts are immediate for the limit geometry: 
\[
 \# \; S_n \ = \ n!, \quad \# \; S_k\times S_{n-k} \ = \ k!(n-k)!, \quad \text{and} \quad \#\; \Sigma(k,n) \ = \ \frac{n!}{k!(n-k)!} \ = \ \binom nk.
\]
The corresponding counts for the geometry of Grassmannians will employ the quantities
\[
 [n]_q \ = \ \sum_{i=0}^{n-1} q^i, \qquad [n]_q! \ = \ \prod_{i=1}^n [i]_q, \qquad \text{and} \qquad \qbinom nk \ = \ \frac{[n]_q!}{[k]_q![n-k]_q!},
\]
which are called the \emph{Gauss number}, the \emph{Gauss factorial} and the \emph{Gauss binomial}, respectively, or sometimes \emph{quantum number}, \emph{quantum factorial} and \emph{quantum binomial} because of their relevance in theoretical physics. Note that we recover the classical quantities in the limit $q\to 1$:
\[
 \lim\limits_{q\to1}\ [n]_q \ = \ n, \qquad \lim\limits_{q\to1}\ [n]_q! \ = \ n!, \qquad \text{and} \qquad \lim\limits_{q\to1}\ \qbinom nk \ = \ \binom nk.
\]
The elements of $G=\GL(n,\F_q)$ correspond to ordered bases of $\F_q^n$. For the first basis vector, we have $q^n-1$ choices, for the second $q^n-q$ choices and so forth. Therefore we have
\[
 \#\GL(n,\F_q) \ = \ \prod_{i=1}^{n} (q^n-q^{i-1}) \ = \ \prod_{i=1}^{n} (q-1) q^{i-1} [i]_q \ = \ (q-1)^n q^{\frac12(n^2-n)} [n]_q!
\]
The elements of $P_k$ decompose into four blocks, consisting of an invertible $k\times k$-matrix, an arbitrary $k\times (n-k)$-matrix, a zero matrix and an invertible $(n-k)\times(n-k)$-matrix. Thus we obtain
\[
 \# P_k \ = \ \#\GL(k,\F_q)\; \#\Mat(k\times n,\F_q)\; \#\GL(n-k,\F_q) \ = \ (q-1)^n q^{\frac12(n^2-n)} [k]_q! [n-k]_q!
\]
Dividing the former two quantities yields
\[
 \# \Gr(k,\F_q^n) \ = \ \frac{\#\GL(n,\F_q)}{\# P_k} \ = \ \frac{[n]_q!}{[k]_q![n-k]_q!} \ = \ \qbinom nk.
\]
Note that the limit $q\to1$ of the cardinalities of $G$ and $P_k$ is $0$ due to the term $(q-1)^n$. But if we resolve this zero, i.e.\ if we divide by $\#T=(q-1)^n$, then we obtain
\[
 \lim\limits_{q\to1}\ \#G/\#T \ = \ \# S_n, \quad \lim\limits_{q\to1}\ \#P_k/\#T \ = \ \# S_k\times S_{n-k}, \quad \lim\limits_{q\to1}\ \#\Gr(k,\F_q^n) \ = \ \# \Sigma(k,n).
\]

Based on these observations, Tits dreamt about the existence of a geometry over a field $\Fun$ with one element that is capable to explain these effects. Later, this idea has been summarized in hypothetical formulas such as
\[
 \GL(n,\Fun) \ = \ S_n \qquad \text{and} \qquad \Gr(n,\Fun) \ = \ \Sigma(k,n).
\]
However, we caution the reader not to take these formulas too literal, for the reasons explained in \cite[Prologue]{L16}.

\begin{ex} \label{ex: Tits geometry of GL_3}
 We illustrate the ideas of Tits in the example of $\GL(3,\F_q)$. In this case, we consider the two Grassmannians $\Gr(1,\F_q^3)$ and $\Gr(2,\F_q^3)$, whose points corresponds to the points and lines in $\P^2(\F_q)$, respectively. The incidence relation consists of pairs of a point $P$ and a line $L$ such that $P\in L$. In the case $q=2$, we calculate
 \[
  \#\Gr(1,\F_2^3) \ = \ \qbinom 31 \ = \ 7 \qquad \text{and} \qquad \#\Gr(2,\F_2^3) \ = \ \qbinom 32 \ = \ 7.
 \]
 Moreover, note that every line contains $q+1=3$ points and that every point is contained in $q+1=3$ lines. The corresponding geometry is illustrated on the left hand side of Figure \ref{figure: geometry of GL3} where the dots correspond to the points in $\P^2(\F_q)$, the circles correspond to the lines in $\P^2(\F_q)$ and an edge between a dot and a circle indicates that the corresponding point is contained in the corresponding line.
 
 Considering the limit $q\to 1$ yields the geometry for $S_3$ that consists of the sets $\Sigma(1,3)=\bigl\{\{1\},\{2\},\{3\}\bigr\}$ and $\Sigma(2,3)=\bigl\{\{1,2\},\{1,3\},\{2,3\}\bigr\}$. Note that every point of this geometry, i.e.\ an one element subset of $\{1,2,3\}$, is contained in $q+1=2$ lines, i.e.\ a $2$-subset of $\{1,2,3\}$. Similarly, every line contains $q+1=2$ points. The corresponding geometry is depicted on the right hand side of Figure \ref{figure: geometry of GL3} .

\newcounter{tikz-counter}

\begin{figure}[ht]
 \[
  \beginpgfgraphicnamed{fig1}
  \begin{tikzpicture}[node distance = 8cm]
   \node (origin) {};
   \foreach \a in {1,...,14}{\draw (\a*360/7+90/7: 2cm) node [draw,circle,inner sep=2pt] (l\a) {};}
   \foreach \a in {1,...,7}{\draw (\a*360/7+270/7: 2cm) node [draw,circle,inner sep=1.8pt,fill=black] (p\a) {};
                              \draw [-] (p\a) -- (l\a);
                              \setcounter{tikz-counter}{\a};
                              \addtocounter{tikz-counter}{1};
                              \draw [-] (p\a) -- (l\arabic{tikz-counter});
                              \addtocounter{tikz-counter}{2};
                              \draw [-] (p\a) -- (l\arabic{tikz-counter});
                             }
   \foreach \a in {1,...,6}{\draw [right of=origin] (\a*360/3+30: 1.2cm) node [draw,circle,inner sep=2pt] (l\a) {};}
   \foreach \a in {1,...,3}{\draw [right of=origin] (\a*360/3+90: 1.2cm) node [draw,circle,inner sep=1.8pt,fill=black] (p\a) {};
                              \draw [-] (p\a) -- (l\a);
                              \setcounter{tikz-counter}{\a};
                              \addtocounter{tikz-counter}{1};
                              \draw [-] (p\a) -- (l\arabic{tikz-counter});
                              \addtocounter{tikz-counter}{2};
                              \draw [-] (p\a) -- (l\arabic{tikz-counter});
                             }
   \draw [->,thick] (3.5,0) -- node[below] {$q\to1$} (5.5,0) ;
  \end{tikzpicture}
  \endpgfgraphicnamed
 \] 
 \caption{The geometry of $\GL(3,\F_q)$ for $q=2$ and its limit $q\to1$}
 \label{figure: geometry of GL3}
\end{figure}
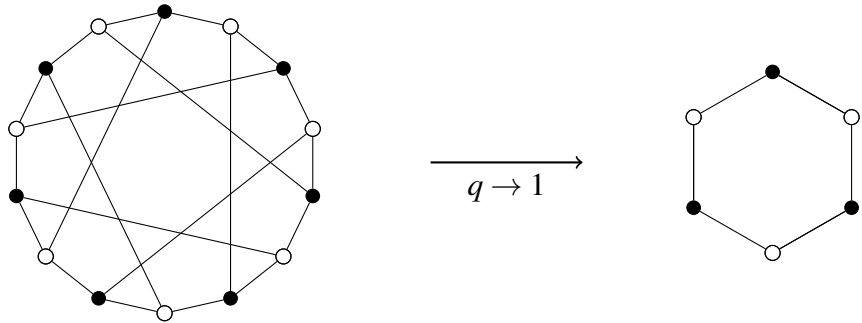
\end{ex}

\subsection{The Riemann hypothesis} \label{subsection: Riemann hypothesis}
One of the most profound problems in mathematics is the Riemann hypothesis. We refrain from alluding to its importance, but refer the reader to one of the numerous overview texts on the topic. A good start is the book \cite{RiemannHypothesis08}.

Arguably, the most influential publication for shaping the area of $\Fun$-geometry was Manin's lecture notes \cite{Manin95}, in which Deninger's programme for a proof of the Riemann hypothesis (\cite{Deninger91}, \cite{Deninger92}, \cite{Deninger92}) and Kurokawa's work on absolute tensor products of zeta functions (\cite{Kurokawa92}) are combined with the idea to realize the integers as a curve over the elusive field $\Fun$ with one element.

In the following, we briefly review the Riemann hypothesis and explain how the proof of its analogue for function fields leads to the desire for algebraic geometry over $\Fun$.


\subsubsection{The Riemann zeta function}
The \emph{Riemann zeta function} is defined by the formulas
\[
 \zeta(s) \ = \ \sum_{n\geq1} \frac{1}{n^s} \ = \ \prod_{p\text{ prime}} \frac{1}{1-p^{-s}},
\]
both of which are expressions that converge absolutely in the half plane $\{s\in\C|\Re s>1\}$. It extends to a meromorphic function on the whole complex plane with simple poles at $0$ and $1$ and it satisfies a certain functional equation. The zeros of the Riemann zeta function are more delicate and the protagonists of the Riemann hypothesis.

Let us begin with explaining the functional equation. The $\Gamma$-function is defined in terms of the formula
\[
 \Gamma(t) \ = \ \int_0^\infty x^{t-1} e^{-x}\; dx,
\]
which converges on the half plane $\{t\in\C | \Re t>0\}$ and can be extended to a meromorphic function on the whole complex plane. It has simple poles at all negative integers $-1,-2,\dotsc$ and no zeros.

The \emph{completed zeta function} is the meromorphic function
\[
 \zeta^\ast(s) \ = \ \Gamma(s/2) \; \pi^{-s/2} \; \zeta(s)
\]
and it satisfies the functional equation
\[
 \zeta^\ast(1-s) \ = \ \zeta^\ast(s).
\]

It is immediate from the definition that $\zeta(s)$ does not have any zero for $\Re s>1$. The functional equation implies that $\zeta(s)$ compensates the poles of $\Gamma(s/2)$ with simple zeros at the negative even integers $-2,-4,\dotsc$, which are called the \emph{trivial zeros} of the zeta function. We see that all other zeros of $\zeta(s)$ lie on the \emph{critical strip $\{s\in\C | 0\leq\Re s\leq 1\}$}. The Riemann hypothesis claims the following.

\begin{RH}
 Every non-trivial zero of $\zeta(s)$ has real part $1/2$.
\end{RH}

\begin{figure}[ht]
 \[
  \beginpgfgraphicnamed{fig2}
  \begin{tikzpicture}[inner sep=0,x=55pt,y=3pt,font=\footnotesize]
   \fill [orange!20] (0,-27) rectangle (1,27);  
   \draw (-4.5,0) -- (1.7,0);
   \draw[dashed] (0,27) -- (0,-27);
   \draw[dashed] (1,27) -- (1,-27);
   \draw[densely dotted] (0.5,27) -- (0.5,-27);
   \draw[fill=white] (0,0) circle (2.5pt);
   \draw[fill=white] (1,0) circle (2.5pt);
   \filldraw (-2,0) circle (2pt);
   \filldraw (-4,0) circle (2pt);
   \filldraw (0.5,14.13) circle (2pt);
   \filldraw (0.5,21.02) circle (2pt);
   \filldraw (0.5,25.01) circle (2pt);
   \filldraw (0.5,-14.13) circle (2pt);
   \filldraw (0.5,-21.02) circle (2pt);
   \filldraw (0.5,-25.01) circle (2pt);
   \draw [decoration={brace},decorate,line width=1pt] (1,-28) -- (0,-28);
   \node (0) at (-2.065,-3.5) {$-2$};
   \node (0) at (-4.065,-3.5) {$-4$};
   \node (0) at (-0.1,-3.5) {$0$};
   \node (0) at (1.1,-3.5) {$1$};
   \node [text width=90pt,align=left] (0) at (1.38,14.13) {$1/2+i\cdot 14.13\dotsc$};
   \node [text width=90pt,align=left] (0) at (1.38,-14.13) {$1/2-i\cdot 14.13\dotsc$};
   \node [text width=90pt,align=left] (0) at (1.38,21.02) {$1/2+i\cdot 21.02\dotsc$};
   \node [text width=90pt,align=left] (0) at (1.38,-21.02) {$1/2-i\cdot 21.02\dotsc$};
   \node [text width=90pt,align=left] (0) at (1.38,25.01) {$1/2+i\cdot 25.01\dotsc$};
   \node [text width=90pt,align=left] (0) at (1.38,-25.01) {$1/2-i\cdot 25.01\dotsc$};
   \node (0) at (0.5,32) {critical line};
   \draw[->] (0.5,30.2) -- (0.5,27.8);
   \node (0) at (0.5,-32) {critical strip};
   \node (0) at (-3,10) {trivial zeros};
   \draw[->] (-3.5,7.5) -- (-3.82,3);
   \draw[->] (-2.5,7.5) -- (-2.18,3);
   \node (0) at (-0.18,1.65) {pole};
   \node (0) at (1.22,1.65) {pole};
   \node [text width=38pt,align=center] (0) at (2.2,0) {\mbox{nontrivial} zeros};
   \draw[->] (1.85,4) -- (1.6,10);
   \draw[->] (1.85,-4) -- (1.6,-10);
  \end{tikzpicture}
  \endpgfgraphicnamed
 \] 
 \caption{The poles and zeros of $\zeta(s)$}
 \label{figure: zeta function}
\end{figure}
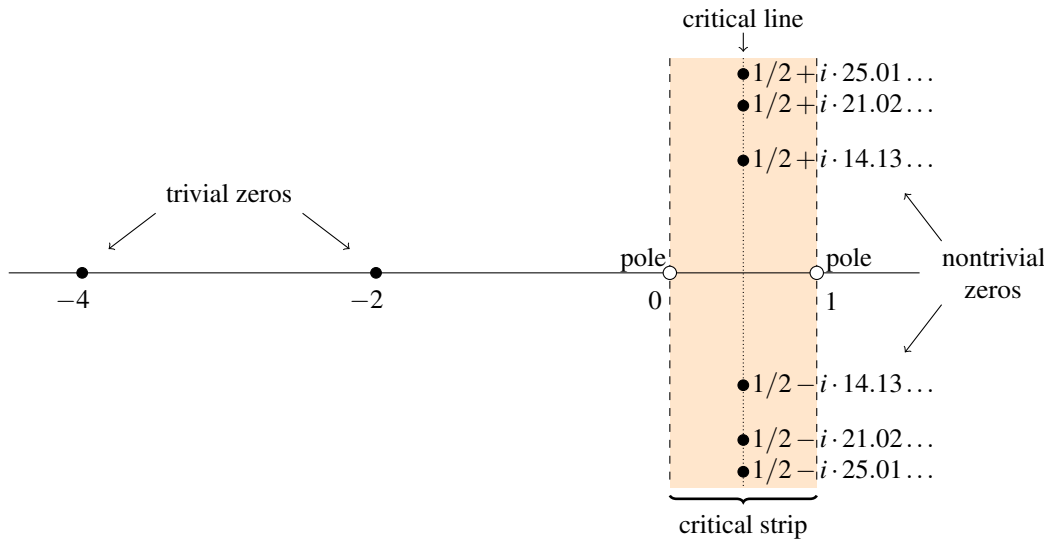

In other words, the Riemann hypothesis states that all zeros of completed zeta function $\zeta^\ast(s)$ lie on the \emph{critical line} $\{s\in\C | \Re s=1/2\}$. See Figure \ref{figure: zeta function} for an illustration of the poles and zeros of $\zeta(s)$.

\subsubsection{Absolute values}
While the Riemann hypothesis is still an open problem, its function field analoguea has been proven by Andr\'e Weil. In order to explain the analogy between the Riemann zeta function $\zeta_\Q(s)=\zeta(s)$ and the zeta function $\zeta_F(s)$ of a function field, we reinterpret the factors of the Euler product $\zeta_\Q(s)=\prod (1-p^{-s})^{-1}$ in terms of absolute values of $\Q$.

An \emph{absolute value} of a field $F$ is a map $v:F\to\R_{\geq0}$ to the non-negative real numbers that satisfies $v(0)=0$, $v(1)=1$, $v(ab)=v(a)v(b)$ and $v(a+b)\leq v(a)+v(b)$ for all $a,b\in F$. 

We illustrate some examples. Every field $F$ has a \emph{trivial absolute value}, which is the absolute value $v_0:F\to\R_{\geq0}$ that maps every nonzero element of $F$ to $1$. For $F=\Q$, the usual absolute value or \emph{archimedean absolute value} $v_\infty=\norm{\ }:\Q\to\R_{\geq0}$ is an absolute value. But there are other, substantially different absolute values for $\Q$. Namely, for every prime number $p$, the \emph{$p$-adic absolute value} is defined as the map $v_p:\Q\to \R_{\geq0}$ with $v_p(p^i\cdot\frac ab)=p^{-i}$ if neither $a$ nor $b$ are divisible by $p$.

An absolute value $v:F\to \R_{\geq0}$ is \emph{nonarchimedean} if it satisfies the \emph{strong triangle inequality} $v(a+b)\leq \max\{v(a),v(b)\}$ for all $a,b\in F$. In fact, most absolute values are nonarchimedean. By a theorem of Ostrowski, the only exceptions come from restricting the usual absolute value of $\C$ to a subfield, and possibly taking powers. If $v$ is nonarchimedean, then the subsets
\[
 \cO_v \ = \ \{ \, a\in F \, | \, v(a) \leq 1 \, \} \qquad \text{and} \qquad \fm_v \ = \ \{ \, a\in F \, | \, v(a) < 1 \, \} 
\]
are a local ring, called the \emph{valuation ring of $v$}, and its maximal ideal, respectively. We define the \emph{residue field of $v$} as $k(v)=\cO_v/\fm_v$. Examples of nonarchimedean places are the trivial absolute value $v_0$ and the $p$-adic absolute value $v_p:\Q\to\R_{\geq0}$. The valuation ring of $v_p$ and its maximal ideal are
\[
 \cO_{v_p} \ = \ \{\, \frac ab\in\Q \, | \, p\text{ does not divide }b\, \} \qquad \text{and} \qquad \fm_{v_p} \ = \ \{\, \frac ab\in\cO_{v_p} \, | \, p\text{ divides }a \, \},
\]
respectively. The residue field $k(v_p)$ is isomorphic to the finite field $\F_p$ with $p$ elements.

Two absolute values $v_1$ and $v_2$ of $F$ are \emph{equivalent} if there is a $t\in\R_{\geq0}$ such that $v_2(a)=v_1(a)^t$ for all $a\in F$. A \emph{place} of $F$ is an equivalence class of nontrivial absolute values. Note that two equivalent nonarchimedean absolute values define the same subsets $\cO_v$ and $\fm_v$, and thus have the same residue field $k(v)$.

By another theorem of Ostrowski, the places of $\Q$ are the \emph{archimedean place}, represented by $v_\infty$, and the \emph{$p$-adic places}, represented by $v_p$, where $p$ ranges through all prime numbers. 

By a slight abuse of notation, we will denote a place represented by an absolute value $v$ by the same symbol $v$. If $v$ is nonarchimedean and $k(v)$ is finite, then we define the \emph{local zeta factor at $v$} as
\[
 \zeta_v(s) \ = \ \frac{1}{1-\#k(v)^s}. 
\]
If we define the zeta factor of the archimedean absolute value $v_\infty$ of $\Q$ as $\zeta_{v_\infty}(s) \ = \ \pi^{-s/2}\Gamma(s/2)$, then the completed zeta function can be expressed as
\[
 \zeta_\Q^\ast(s) \ = \ \prod_{\text{places }v\text{ of }\Q} \zeta_v(s)
\]
in the region of convergence, i.e.\ for $\Re s>1$. While the different shape of the factor $\zeta_\infty$ at infinity has been the cause for much musing, the nonarchimedean factors have a direct analogue in the function field setting.

\subsubsection{Zeta functions for function fields}
Let $F$ be a global function field, i.e.\ a finite field extension of a rational function field $\F_q(T)$ over a finite field $\F_q$ with $q$ elements. For later reference, we assume that $q$ is maximal with the property that $F$ contains $\F_q(T)$ as a subfield. Every nontrivial absolute value of $F$ is nonarchimedean and the residue field is a finite field extension of $\F_q$ and therefore finite. Thus the local zeta factors $\zeta_v(s)=(1-\#k(v)^{-s})^{-1}$ make sense, and we define the zeta function of $F$ as
\[
 \zeta_F(s) \ = \ \prod_{\text{places }v\text{ of }F} \zeta_v(s),
\]
which is an expression that has analogous properties to the completed Riemann zeta function: it converges for $\Re s>1$ and has a meromorphic continuation to all $s$ in $\C$. It satisfies a functional equation of the form
\[
 \zeta_F(1-s) \ = \ \pm q^{(2g-2)(1-s)} \zeta_F(s)
\]
where $g$ is the genus of $F$, a number which plays an analogue role as the genus of a Riemann surface. As explained below, the field $F$ occurs indeed as the function field of a certain curve, but explaining the definition of the genus would lead us too far astray. 

A particular property for the function field setting is that the zeta function can be expressed in terms of a rational function. Namely, $\zeta_F(s)=Z_F(q^{-s})$ for a function $Z_F(T)$ in $T$ of the form
\[
 Z_F(T) \ = \ \frac{P(X)}{(1-T)(1-qT)}.
\]
This implies that $\zeta_F(s)$ is periodic modulo $\frac{2\pi i}{\ln q}$, and it has simple poles in all complex numbers of the form $k\frac{2\pi i}{\ln q}$ and $1+k\frac{2\pi i}{\ln q}$ with $k\in\Z$; cf.\ the illustration in Figure \ref{figure: zeta function for function fields} below.

\begin{ex}
 To explain the analogy to $\Q$ in more detail below, we consider the example of a rational function field $F=\F_q(T)$. Then $F$ is the field of fraction of the polynomial ring $\F_q[T]$. 
 
 The analogue of the $p$-adic absolute value of $\Q$ is the $f$-adic absolute value $v_f:\F_q(T)\to\R_{\geq0}$ where $f$ is an irreducible polynomial in $\F_q[T]$, i.e.\ $f$ is a polynomial of positive degree $d\geq1$ and does not equal the product of two polynomials of positive degree. Note that we can write every nonzero element of $\F_q(T)$ in the form $f^i\frac gh$ where $i\in\Z$ and $g$ and $h$ are polynomials that are not multiples of $f$ by another polynomial. Then the $f$-adic valuation $v_f:\F_q(T)\to\R_{\geq0}$ is defined by the formula
 \[
  v_f(f^i\frac gh) \ = \ q^{-di}.
 \]
 Note that, similar to the $p$-adic valuation, $v_f$ is nonarchimedean and we have
 \[
  \cO_{v_f} \ = \ \{\, \frac gh \in\F_q(T) \, | \, f\text{ does not divide }h\, \} \quad \text{and} \quad \fm_{v_f} \ = \ \{\, \frac gh\in\cO_{v_f} \, | \, f\text{ divides }g \, \}.
 \]
 The residue field $k(v_f)=\cO_{v_f}/\fm_{v_f}$ is isomorphic to $\F_{q^d}$, the unique degree $d$-extension field of $\F_q$.
 
 Up to equivalence, the $f$-adic absolute values $v_f$ represent all places of $\F_q(T)$, with the exception of the place at infinity. This latter place is represented by the absolute value $v_\infty$, which is defined by
 \[
  v_\infty(\frac gh) \ = \ q^{\deg h-\deg g}
 \]
 for nonzero polynomials $g$ and $h$ of respective degrees $\deg g$ and $\deg h$. The absolute value $v_\infty$ is also nonarchimedean, in stark contrast to the situation for $\Q$ where the place at infinity is archimedean. We have
 \[
  \cO_{v_\infty} \ = \ \{\, \frac gh \in\F_z(T) \, | \, \deg g\leq\deg h\, \} \quad \text{and} \quad \fm_{v_\infty} \ = \ \{\, \frac gh\in\cO_{v_\infty} \, | \, \deg g < \deg h \, \}.
 \]
 The residue field of $v_\infty$ is $k(v_\infty)=\F_q$.
 
 To conclude this example, we note that it is not too hard to show that the zeta function of $F=\F_q(T)$ has the following form:
 \[
  \zeta_F(s) \ = \ \prod_{\text{places }v\text{ of }F} \frac{1}{1-\#k(v)^{-s}} \ = \ \frac{1}{(1-q^{-s})(1-q^{1-s})}.
 \]
 In particular, $\zeta_F(s)$ does not have any zero at all in this case.
\end{ex}

\subsubsection{From the Hasse-Weil theorem to \texorpdfstring{$\Fun$}{F1}}
The analogue of the Riemann hypothesis for $\zeta_F(s)$ has been proven by Hasse (\cite{Hasse36}) in the case of elliptic function fields and by Weil (\cite{Weil48}) for all function fields. It is also known as the \emph{Hasse-Weil theorem}.

\begin{thm}
 Every zero of $\zeta_F(s)$ has real part $1/2$.
\end{thm}


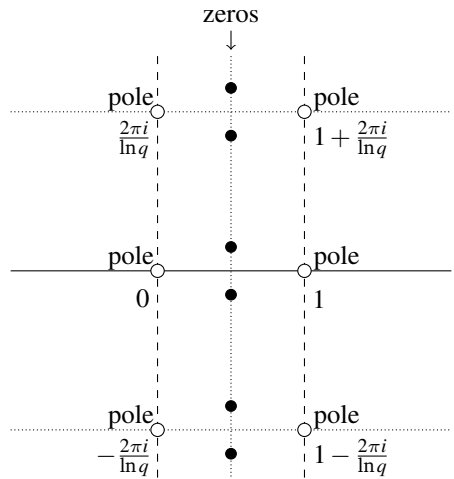
\begin{figure}[ht]
 \[
  \beginpgfgraphicnamed{fig3}
  \begin{tikzpicture}[inner sep=0,x=55pt,y=3pt,font=\footnotesize]
   \draw (-1,0) -- (2,0);
   \draw[densely dotted] (-1,20) -- (2,20);
   \draw[densely dotted] (-1,-20) -- (2,-20);
   \draw[dashed] (0,27) -- (0,-27);
   \draw[dashed] (1,27) -- (1,-27);
   \draw[densely dotted] (0.5,27) -- (0.5,-27);
   \draw[fill=white] (0,0) circle (2.5pt);
   \draw[fill=white] (1,0) circle (2.5pt);
   \draw[fill=white] (0,20) circle (2.5pt);
   \draw[fill=white] (1,20) circle (2.5pt);
   \draw[fill=white] (0,-20) circle (2.5pt);
   \draw[fill=white] (1,-20) circle (2.5pt);
   \filldraw (0.5,3) circle (2pt);
   \filldraw (0.5,17) circle (2pt);
   \filldraw (0.5,23) circle (2pt);
   \filldraw (0.5,-3) circle (2pt);
   \filldraw (0.5,-17) circle (2pt);
   \filldraw (0.5,-23) circle (2pt);
   \node (0) at (-0.1,-3.5) {$0$};
   \node (0) at (1.1,-3.5) {$1$};
   \node (0) at (-0.16,16.5) {$\frac{2\pi i}{\ln q}$};
   \node (0) at (1.32,16.5) {$1+\frac{2\pi i}{\ln q}$};
   \node (0) at (-0.23,-23.5) {$-\frac{2\pi i}{\ln q}$};
   \node (0) at (1.32,-23.5) {$1-\frac{2\pi i}{\ln q}$};
   \node (0) at (0.5,32) {zeros};
   \draw[->] (0.5,30.2) -- (0.5,27.8);
   \node (0) at (-0.18,1.65) {pole};
   \node (0) at (1.22,1.65) {pole};
   \node (0) at (-0.18,21.65) {pole};
   \node (0) at (1.22,21.65) {pole};
   \node (0) at (-0.18,-18.35) {pole};
   \node (0) at (1.22,-18.35) {pole};
  \end{tikzpicture}
  \endpgfgraphicnamed
 \] 
 \caption{The poles and zeros of $\zeta_F(s)$}
 \label{figure: zeta function for function fields}
\end{figure}

\begin{ex}
 To give a few examples, the case of genus $g=0$ is the case of a rational function field $F=\F_q(T)$ where the curve $C$ is a projective line over $\F_q$. As we have seen before, the Hasse-Weil theorem is trivial in this case since $\zeta_F(s)$ does not have any zeroes. The case of genus $g=1$ is the case of an elliptic curve, which has been treated by Hasse. In this case, the Riemann zeta function has two zeros modulo $\frac{2\pi i}{\ln q}$ and the Riemann hypothesis for $F$ is equivalent to the estimate $q-2\sqrt{q}\leq h\leq q+2\sqrt{q}$ where $h$ is the number of places $v$ of $F$ with residue field $k(v)$ equal to $\F_q$.
\end{ex}

In the following, we will explain a few of the key ingredients in the proof of the Hasse-Weil theorem and make clear how this leads to the postulation of $\Fun$.

As mentioned before, $F$ is the function field of a curve $C$ over $\F_q$. The points of $C$ are, by definition, the places of $F$ where we neglect the ``generic point'' in our description. The curve fibres over its ``base point'' $\Spec \F_q$. 

We can embed the curve $C$ diagonally into the fibre product $C\times_{\Spec\F_q}C$. The Riemann hypothesis for $\zeta_F(s)$ can be tied to an estimate for the number of intersection points of the diagonally embedded curve $C$ with its ``Frobenius twist'' inside the surface $C\times_{\Spec\F_q}C$. This estimate can be established by an explicit calculation.

The analogies between number fields and function fields lead to the hope that one can mimic these methods for $\Q$ and approach the Riemann hypothesis. Grothendieck's theory of schemes provides a satisfying framework to view the collection of all nonarchimedean places $v_p$ of $\Q$ as a curve, namely, as the \emph{arithmetic curve} $\Spec\Z$, the spectrum of $\Z$. In analogy to the function field setting, we would like to include the archimedean place $v_\infty$ in a hypothetical \emph{completion $\overZ$ of $\Spec\Z$} and we would like to have a base field for $\overZ$, namely $\Fun$, the field with one element.

In such a theory, we should have a base extension functor $-\otimes_\Fun\Z$ from $\Fun$-schemes to usual schemes and we should be able to define the \emph{arithmetic surface} $\overZ\times_\Fun\overZ$.


\section{Birth hour: \texorpdfstring{$\Fun$}{F1}-varieties and their siblings}
\label{section: birth hour}

The first person courageous enough to propose a notion of an $\Fun$-variety was Cristophe Soul\'e who presented a first attempt \cite{Soule99} at the Mathematische Arbeitstagung of the MPI in 1999. He published a refinement \cite{Soule04} of his approach in 2004. 

Soon after, ideas were blooming and the same decade saw more than a dozen different definitions for $\Fun$-geometry. In the following, we give a brief and incomplete overview of such approaches. For more details and references, confer \cite{LopezPena-L11b,L16}.

Soul\'e's approach underwent a number of further variations by himself in \cite{Soule11} and by Connes and Consani in \cite{Connes-Consani11a,Connes-Consani10a}. These approaches are closely related to the notion of a torified scheme as introduced by L\'opez Pe\~na and the author in \cite{LopezPena-L11a}.

One of the most important approaches is Deitmar's variation \cite{Deitmar05} of Kato fans (\cite{Kato94}), which he calls $\Fun$-schemes. It is the most minimalistic approach to $\Fun$-geometry since it is included in every other theory about $\Fun$. In this sense, Deitmar's theory constitutes the very core of $\Fun$-geometry. Subsequently, this theory has found various applications under the names of Deitmar schemes, monoid schemes or monoidal schemes.

To\"en and Vaqui\'e generalize in \cite{Toen-Vaquie09} the functorial viewpoint on scheme theory to any category $\cC$ that looks sufficiently like a categories of modules over rings. Durov gains a notion of $\Fun$-varieties in his unpublished approach \cite{Durov07} to Arakelov theory. Borger defines in \cite{Borger09} an $\Fun$-variety as a scheme together with a lift of the Frobenius automorphisms for every prime number $p$. Haran has developed various approaches in \cite{Haran07,Haran10,Haran17}. 

Lescot develops in \cite{Lescot11,Lescot12} an algebraic geometry over idempotent semirings, which he dubs $\Fun$-geometry. Connes and Consani extend Lescot's viewpoint to the context of Krasner's hyperrings, which they promote as a geometry over $\Fun$ in \cite{Connes-Consani10b,Connes-Consani11b}. This point of view has been extended by Jun in \cite{Jun15}.

Berkovich introduces a theory of congruence schemes for monoids (\cite{Berkovich11}), which can be seen as an enrichment of Deitmar's $\Fun$-geometry. Deitmar modifies this approach in \cite{Deitmar11a}.

The author develops in \cite{L12a,L12d,LopezPena-L12,L14b,L17} (partly in collaboration with L\'opez Pe\~na) a notion of $\Fun$-geometry, based on the notion of a so-called blueprint.

We will examine Deitmar's and the author's approach to $\Fun$-geometry in more detail in the following sections.


\subsection{Monoid schemes}\label{subsection: monoid schemes}
In this section, we review a slight modification of Deitmar's approach to $\Fun$-schemes in \cite{Deitmar05}, which can be seen as the very core of $\Fun$-geometry. It realizes the motto ``non-additive geometry'' literally and it forms a subclass of every other approach to varieties over $\Fun$, up to some finiteness conditions in certain cases. We make the definitions approachable to the non-expert, but warn the inexperienced reader that a motivation for scheme theory lies outside the scope of this overview paper.

\subsubsection{Monoids}
The underlying algebraic objects in \cite{Deitmar05} are monoids, i.e.\ semigroups with an identity element. For the purpose of $\Fun$-geometry, it has been proven useful to consider \emph{monoids with zero}, which are monoids $A$ together with an absorbing element $0$, i.e. $0\cdot a=0$ for all $a\in A$ where we write the monoid multiplicatively. In this exposition, we will consider the variation of monoid schemes for monoids with zero.

In the following, we will agree that all of our monoids are commutative and with zero. A \emph{monoid morphism} is a map $f:A_1\to A_2$ between monoids $A_1$ and $A_2$ with $f(0)=0$, $f(1)=1$ and $f(ab)=f(a)f(b)$ for all $a,b\in A_1$.

There is a base extension functor $-\otimes_\Fun\Z$ that sends a monoid to a ring. Namely, given a monoid $A$, we define 
\[
 A\otimes_\Fun\Z \ = \ \Z[A] \, / \, \gen{0_A}
\]
where $\Z[A]$ is the semigroup ring $\Z[A]=\{\sum n_a a|n_a\in \Z,\text{ almost all }0\}$ of finite $\Z$-linear combinations of elements of $A$ and where $\gen{0_A}$ is the ideal generated by the zero $0_A$ of $A$. In other words, we identify the zero of the monoid $A$ with the zero of the ring $\Z[A]$. A monoid morphism $f:A_1\to A_2$ can be extended by linearity to a ring homomorphism $A_1\otimes_\Fun\Z\to A_2\otimes_\Fun\Z$ between the associated rings, which defines the base extension functor $-\otimes_\Fun\Z$ from monoids to rings.

\begin{ex}
 We provide some first examples of monoids. The trivial monoid is the monoid with a single element $0=1$. Its base extension $\{0\}\otimes_\Fun\Z$ is the trivial ring.
 
 The smallest nontrivial monoid consists solemnly of two elements $0$ and $1$, and it is this monoid that we call $\Fun$. Its base extension $\Fun\otimes_\Fun\Z$ is $\Z$.
 
 The free monoid generated by a number of indeterminants $T_1,\dotsc,T_n$ consists of all monomials $T_1^{e_1}\dotsc T_n^{e_n}$ in $T_1,\dotsc,T_n$ together with a distinct element $0$. We denote the free monoid in $T_1,\dotsc,T_n$ by $\Fun[T_1,\dotsc,T_n]$. Its base extension $\Fun[T_1,\dotsc,T_n]\otimes_\Fun\Z$ is the polynomial ring $\Z[T_1,\dotsc,T_n]$.
 
 Note that every ring is a monoid if we omit its addition.
\end{ex}


\subsubsection{The spectrum}
Let $A$ be a monoid. An \emph{ideal of $A$} is a multiplicative subset $I$ of $A$ such that $0\in I$ and $IA=I$. A \emph{prime ideal of $A$} is an ideal $\fp$ of $A$ such that its complement $S=A-\fp$ is a \emph{multiplicative subset}, i.e.\ it contains $1$ and is closed under multiplication. The \emph{spectrum $\Spec A$ of $A$} is the set of all prime ideals of $A$ together with a topology and a structure sheaf, which we will describe below.

The topology of $\Spec A$ is generated by the \emph{principal open subsets}
\[
 U_h \ = \ \{ \, \fp\in\Spec A \, | \, h\notin \fp \, \}
\]
where $h$ ranges through all elements of $A$. Note that $U_1=\Spec A$, that $U_0=\emptyset$ and that $U_{gh}=U_g\cap U_h$. Thus every open subset of $\Spec A$ is a union of principal open subsets.

\begin{ex}\label{ex: monoid spectra}
 We describe the topological spaces of some spectra of monoids. We begin with some general observations that are helpful in the calculation of the prime ideals. The \emph{unit group} $A^\times$ of a monoid $A$ is the set of invertible elements and forms an abelian group. If $a\in A^\times$ and $\fp$ is a prime ideal, then $a\notin \fp$. Note that $A-A^\times$ is always a prime ideal and that it contains every other prime ideal. 
 
 As a consequence, a monoid $A$ with $A^\times=A-\{0\}$ has a single prime ideal, which is $\{0\}$. In particular, $\Spec\Fun$ consists of a single point.
 
 The ideal $I=(J)$ generated by a subset $J$ of $A$ is, by definition, the smallest ideal containing $J$. If $A$ is generated by elements $a_1,\dotsc,a_n$ as a monoid, i.e.\ every element of $A$ is a finite product of these elements or $0$, then every prime ideal of $A$ is generated by a subset of the generators. 
 
 The free monoid $\Fun[T_1,\dotsc,T_n]$ is generated by $T_1,\dotsc,T_n$. In this case, every subset of $\{T_1,\dotsc,T_n\}$ generates a prime ideal. We illustrate the spectra of $\Fun[T]$ and $\Fun[T_1,T_2]$ in Figure \ref{figure: affine spaces}. The labelled dots stay for the corresponding prime ideals, and a line between two different prime ideals indicates that the prime ideal on the bottom end is contained in the prime ideal on the top end. 
 
 Note that all our examples are topological spaces with finitely many points. In this case, the open subsets are precisely those that are closed from below with respect to the inclusion relation, and the closed subsets are those that are closed from above.
 
 \begin{figure}[ht]
  \[
   \beginpgfgraphicnamed{fig4}
   \begin{tikzpicture}[inner sep=0,x=30pt,y=20pt,font=\footnotesize]
    \filldraw (0,1) circle (3pt);
    \filldraw (0,3) circle (3pt);
    \draw (0,1) -- (0,3);
    \node at (0.5,3) {$(T)$};
    \node at (0.5,1) {$\{0\}$};
    \filldraw (5,2) circle (3pt);
    \filldraw (6,0) circle (3pt);
    \filldraw (6,4) circle (3pt);
    \filldraw (7,2) circle (3pt);
    \draw (5,2) -- (6,0) -- (7,2) -- (6,4) -- (5,2);
    \node at (6.6,0) {$\{0\}$};
    \node at (6.8,4) {$(T_1,T_2)$};
    \node at (7.6,2) {$(T_2)$};
    \node at (4.4,2) {$(T_1)$};
   \end{tikzpicture}
   \endpgfgraphicnamed
  \] 
  \caption{The spectra of $\Fun[T]$ and $\Fun[T_1,T_2]$}
  \label{figure: affine spaces}
\end{figure}
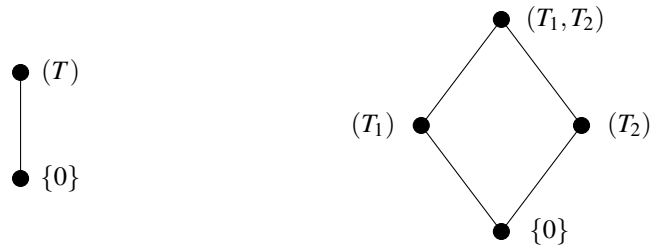
\end{ex}

\subsubsection{The structure sheaf}
It is somewhat more difficult to describe the structure sheaf $\cO_X$ of $X=\Spec A$. It is an association that sends an open subset $U$ of $X$ to a monoid $\cO_X(U)$ and an inclusion $V\subset U$ to a monoid morphism $\res_{U,V}:\cO_X(U)\to\cO_X(V)$. We will restrict ourselves to an explicit description of its values on principal open subsets, which are given in terms of localizations. 

Namely, let $A$ be a monoid and $S$ a multiplicative subset. The \emph{localization of $A$ at $S$} is the monoid
\[
 S^{-1}A \ = \ S\times A / \sim
\]
where $\sim$ is the equivalence relation on the Cartesian product $S\times A$ given by $(s,a)\sim (s',a')$ if and only if there exists a $t\in S$ such that $tsa'=ts'a$. If we denote the equivalence class of $(s,a)$ by $\frac as$, then the multiplication of $S^{-1}A$ is given by the rule $\frac as\cdot \frac bt=\frac{ab}{st}$. The identity element of $S^{-1}A$ is $\frac 11$ and its zero is $\frac 01$. Note that the localization comes with a monoid morphism $A\to S^{-1}A$, defined by $a\mapsto \frac a1$, which sends every $s\in S$ to an invertible element in $S^{-1}A$. To wit, the inverse of $\frac s1$ is $\frac 1s$.

Cases of particular importance are the localizations $A[h^{-1}]=S^{-1}A$ of $A$ in an element $h$ where $S=\{h^i\}_{i\geq0}$ and the localization $A_\fp=S^{-1}A$ of $A$ at a prime ideal $\fp$ where $S=A-\fp$.

\begin{ex}
 Let $A=\Fun[T]$ be the free monoid in one indeterminate $T$. We denote the localization $A[T^{-1}]$ of $\Fun[T]$ in $T$ by $\Fun[T^{\pm 1}]$, which consists in $0$ and all integer powers $T^i$ of $T$. The base extension $\Fun[T^{\pm1}]\otimes_\Fun\Z$ is the ring $\Z[T^{\pm1}]$ of Laurent polynomials over $\Z$. 
\end{ex}

We are prepared to describe the structure sheaf for principal open subsets. Namely, we define $\cO_X(U_h)=A[h^{-1}]$. It follows easily from the definitions that $U_g\subset U_h$ if and only if $g$ is divisible by $h$, i.e.\ we have $g=fh$ for some $f\in A$. This observation lets us define the monoid morphism
\[
 \begin{array}{cccc}
  \res_{U_h,U_g}: & \cO_X(U_h)=A[h^{-1}] & \longrightarrow & A[g^{-1}] = \cO_X(U_g).\\
                  & \frac{a}{h^i}        & \longmapsto     & \frac{af^i}{g^i}
 \end{array}
\]
For example, if $h=1$, then $U_h=\Spec A$ and $\cO_X(\Spec A)=A$. In this case, the restriction map $\res_{U_h,U_g}:A\to A[g^{-1}]$ is nothing else than the aforementioned map that sends $a\to \frac a1$.


The structure sheaf $\cO_X$ of $X=\Spec A$ derives from our description on its values on principal open subsets due to the fact that this definition extends uniquely to a sheaf of monoids\footnote{In order to avoid a digression into technicalities, we do not introduce sheaves. The reader can safely omit all details concerning sheaves.}
on $\Spec A$. Thus $\Spec A$ is a \emph{monoidal space}, i.e.\ a topological space together with a sheaf in monoids. 

\subsubsection{Monoid schemes}
A \emph{monoid scheme}, or \emph{$\Fun$-scheme} in the terminology of \cite{Deitmar05}, is a monoidal space that has a covering by open subsets that are isomorphic to the spectra of monoids. We say that a monoid scheme is \emph{affine} if it is isomorphic to the spectrum of a monoid.

The base extension functor $-\otimes_\Fun\Z$ extends to a functor from monoid schemes to usual schemes in terms of open coverings and the definition 
\[
 (\Spec A)\otimes_\Fun\Z \ = \ \Spec \bigl(A\otimes_\Fun\Z\bigr).
\]

\begin{ex}
 In analogy to the affine space $\A^n_\Z=\Spec\Z[T_1,\dotsc,T_n]$ over $\Z$, we define the affine space $\A^n_\Fun$ over $\Fun$ as $\Spec\Fun[T_1,\dotsc,T_n]$. We obtain $\A^n_\Fun\otimes_\Fun\Z=\A^n_\Z$. 
 
 Similarly, we define the multiplicative group scheme $\G_{m,\Fun}$ over $\Fun$ as $\Spec\Fun[T^{\pm1}]$ in analogy to the multiplicative group scheme $\G_{m,\Z}=\Spec\Z[T^{\pm1}]$ over $\Z$. We obtain $\G_{m,\Fun}\otimes_\Fun\Z=\G_{m,\Z}$. For details about the group law of $\G_{m,\Fun}$, see section \ref{subsection: algebraic groups over F1}.
\end{ex}

This view on $\Fun$-geometry is very appealing thanks to its simple and clean approach. However, it is often too limited for applications. Namely, the only varieties that are base extensions of monoid schemes are toric varieties. In more detail, Deitmar proves the following in \cite{Deitmar08}.


\begin{thm}\label{thm: toric variety}
 Let $X$ be a monoid scheme such that $X\otimes_\Fun\Z$ is a connected, separated and flat scheme of finite type. Then $X\otimes_\Fun\Z$ is a toric variety over $\Z$.
\end{thm}

\begin{ex}\label{ex: monoid schemes}
 As a first example of a monoid scheme that is not affine, we consider the projective line $P_\Fun^1$ over $\Fun$. It consists of two closed points $x_0$ and $x_1$ and a so-called generic point $\eta$, which is contained in every non-empty open subset. It can be covered by the two open subsets $U_0=\{x_0,\eta\}$ and $U_1=\{x_1,\eta\}$, which can be identified with affine lines $\Spec\Fun[T_i]$ over $\Fun$ where $i\in\{0,1\}$. Their intersection is $U_{01}=\{\eta\}$, which coincides with $\G_{m,\Fun}=\Spec\bigl(\Fun[T_0,T_1]/(T_0T_1=1)\bigr)$. The restriction maps from $U_i$ to $U_{01}$ are the obvious inclusions
 \[
  \res_{U_i,U_{01}}: \Fun[T_i] \quad \longrightarrow \quad \Fun[T_0,T_1]/(T_0T_1=1).
 \]
 
 The projective line $\P^1_\Fun$ is illustrated on the left hand side of Figure \ref{figure: projective spaces} where we label the different points with homogeneous coordinates with coefficients in $\Fun$. However, a coefficient $1$ should be read as a generic value different from zero. This explains why $[1:1]$ is the generic point.

 In a similar vain, it is possible to define projective spaces or any toric variety as a monoid scheme. We illustrate the projective surface $\P^2_\Fun$ over $\Fun$ on the right hand side of Figure \ref{figure: projective spaces}. 

 \begin{figure}[ht]
  \[
   \beginpgfgraphicnamed{fig5}
   \begin{tikzpicture}[inner sep=0,x=30pt,y=20pt,font=\footnotesize,back line/.style={},cross line/.style={preaction={draw=white, -,line width=6pt}}]
    \filldraw (1,1) circle (3pt);
    \filldraw (0,3) circle (3pt);
    \filldraw (2,3) circle (3pt);
    \draw (0,3) -- (1,1) -- (2,3);
    \draw (6,4) -- (6,2) -- (8,4) -- (10,2) -- (10,4);
    \draw [cross line] (6,4) -- (8,2) -- (10,4);
    \draw (8,0) -- (6,2);
    \draw (8,0) -- (8,2);
    \draw (8,0) -- (10,2);
    \filldraw (8,0) circle (3pt);
    \filldraw (6,2) circle (3pt);
    \filldraw (8,2) circle (3pt);
    \filldraw (10,2) circle (3pt);
    \filldraw (6,4) circle (3pt);
    \filldraw (8,4) circle (3pt);
    \filldraw (10,4) circle (3pt);
    \node at (1.6,1) {$[1:1]$};
    \node at (-0.6,3) {$[1:0]$};
    \node at (2.6,3) {$[0:1]$};
    \node at (8.9,0) {$[1:1:1]$};
    \node at (5.2,2) {$[1:1:0]$};
    \node at (8.9,2) {$[1:0:1]$};
    \node at (10.8,2) {$[0:1:1]$};
    \node at (5.2,4) {$[1:0:0]$};
    \node at (8.9,4) {$[0:1:0]$};
    \node at (10.8,4) {$[0:0:1]$};
   \end{tikzpicture}
   \endpgfgraphicnamed
  \] 
  \caption{The projective line $\P^1_\Fun$ and the projective plane $\P^2_\Fun$ over $\Fun$}
  \label{figure: projective spaces}
 \end{figure}
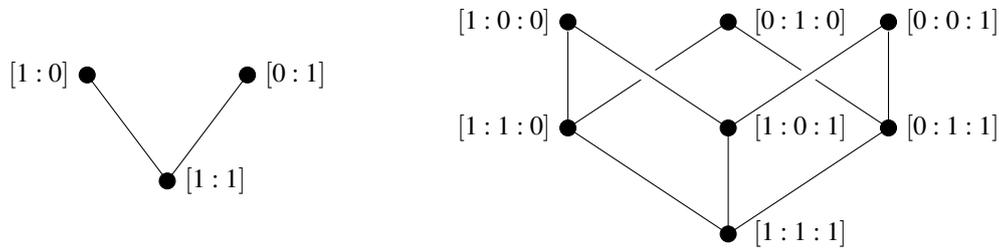
\end{ex}

\begin{rem}
 We see in these examples already a first relation to the incidence geometries considered by Tits. The homogeneous coordinates $[x_1:\dotsb:x_n]$ of a point in $\P^{n-1}_\Fun$ defines a the subset $I=\{i|x_i\neq0\}$ of $\{1,\dotsc,n\}$ and thus a point in $\Sigma(k,n)$ if $k=\# I$. If we consider the action of $S_n$ on $P^{n-1}_\Fun$ that permutes the coordinates of points, then we obtain the limit geometry of $S_n$ as considered in section \ref{subsection: incidence geometry}.

 \begin{figure}[b]
  \[
   \beginpgfgraphicnamed{fig6}
   \begin{tikzpicture}[inner sep=0,x=30pt,y=20pt,font=\footnotesize,back line/.style={},cross line/.style={preaction={draw=white, -,line width=6pt}}]
    \draw (6,2) -- (8,4) -- (10,2);There has 
    \draw [cross line] (6,4) -- (8,2) -- (10,4);
    \draw (6,4) -- (6,2);
    \draw (10,2) -- (10,4);
    \filldraw[fill=white] (6,2) circle (3pt);
    \filldraw[fill=white] (8,2) circle (3pt);
    \filldraw[fill=white] (10,2) circle (3pt);
    \filldraw (6,4) circle (2.5pt);
    \filldraw (8,4) circle (2.5pt);
    \filldraw (10,4) circle (2.5pt);
    \draw[thick,<->] (11.8,3) -- (13.2,3);
    \foreach \a in {1,...,6}{\draw (16,3)+(\a*360/3+30: 1.2cm) node [draw,circle,inner sep=2pt] (l\a) {};}
    \foreach \a in {1,...,3}{\draw (16,3)+(\a*360/3+90: 1.2cm) node [draw,circle,inner sep=1.8pt,fill=black] (p\a) {};
                              \draw [-] (p\a) -- (l\a);
                              \setcounter{tikz-counter}{\a};
                              \addtocounter{tikz-counter}{1};
                              \draw [-] (p\a) -- (l\arabic{tikz-counter});
                              \addtocounter{tikz-counter}{2};
                              \draw [-] (p\a) -- (l\arabic{tikz-counter});
                             }
    \end{tikzpicture}
   \endpgfgraphicnamed
  \] 
  \caption{The relation of $\P^2_\Fun$ to the Tits geometry of $S_3$}
  \label{figure: projective surface and Tits geometry}
 \end{figure}
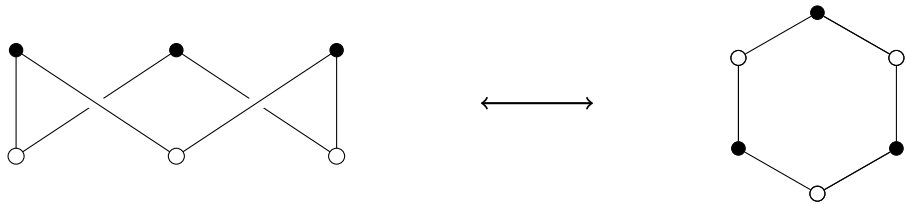

 We illustrate this relation in Figure \ref{figure: projective surface and Tits geometry}. If we remove the generic point $[1:1:1]$ from $\P^2_\Fun$, which corresponds to the full subset of $\{1,2,3\}$, then we are left with the space illustrated on the left hand side of Figure \ref{figure: projective surface and Tits geometry} where we indicate the points corresponding to $2$-subsets by circles. This is the same as the limit geometry from Figure \ref{figure: geometry of GL3}, as repeated on the right hand side of Figure \ref{figure: projective surface and Tits geometry}.
\end{rem}


\subsection{Blueprints and blue schemes}\label{subsection: blueprints}

Since semi-simple algebraic groups are not toric varieties, Theorem \ref{thm: toric variety} testifies that monoid schemes are not sufficient to realize Tits' dream of algebraic groups over $\Fun$. This led to the refinement of monoid schemes in terms of blueprints.

As a leading example for the exposition in this section, we consider the special linear group $\SL(2)$. As a scheme over the integers, we have 
\[
 \SL(2)_\Z \ = \ \Spec\,\bigl(\,\Z[T_1,T_2,T_3,T_4] \, \bigl/ \, (T_1T_4-T_2T_3-1)\, \bigr),
\]
which encodes all $2\times 2$-matrices with coefficients $T_1,\dotsc,T_4$ with determinant $T_1T_4-T_2T_3$ equal to $1$. We would like to make sense of the formula
\[
 \SL(2)_\Fun \ = \ \Spec\,\bigl(\,\Fun[T_1,T_2,T_3,T_4] \, \bigl/ \, (T_1T_4=T_2T_3+1)\, \bigr).
\]

\subsubsection{Semirings}
In order to understand the concept of a blueprint, we have to review some facts about semirings. Although the axioms of a semiring are very similar to those of a ring---we just omit the axiom about additive inverses---, we will see that the theory of semirings shows certain effects that do not occur in ring theory.

A \emph{commutative semiring with $0$ and $1$}, or for short a \emph{semiring}, is a set $R$ together with an addition $+$, a multiplication $\cdot$ and constants $0$ and $1$ such that $(R,+)$ is a commutative semigroup with neutral element $0$, such that $(R,\cdot)$ is a commutative semigroup with neutral element $1$ and with zero $0$ and such that $a(b+c)=ab+ac$ for all $a,b,c\in R$. A \emph{semiring homomorphism} is an additive and multiplicative map between semirings that maps $0$ to $0$ and $1$ to $1$.

\begin{ex}\label{ex: semirings}
 Commutative rings with $1$ are semirings. Other examples are the natural numbers $\N$ and the non-negative real numbers $\R_{\geq0}$ together with the usual addition and multiplication. 
 
 Given a semiring $R$, we can form the polynomial semiring $R[T_1,\dotsc,T_n]$ in $n$ indeterminants $T_1,\dotsc,T_n$, which consists in all polynomials $\sum a_e T^e$ where $e=(e_1,\dotsc,e_n)$ is a multi-index and $T^e=T_1^{e_1}\dotsb T_n^{e_n}$. Addition and multiplication of $R[T_1,\dotsc,T_n]$ are defined in the same way as for polynomial rings.
 
 Given a multiplicative monoid $A$ with $0$, then we define the monoid semiring as the set $\N[A]=\bigl\{\sum a_i\bigl|a_i\in A-\{0\}\bigl\}$ of finite formal sums of nonzero elements in $A$ whose product is inherited by the product of $A$. Note that if $ab=0$ in $A$, then $1a\cdot 1b$ equals the neutral element for addition of $\N[A]$, which is the empty sum. In other words, the zero of $A$ is identified with the zero of $\N[A]$.
 
 An example of a more exotic, or \emph{tropical}, nature is the Boolean semifield $\B=\{0,1\}$ where $1+1=1$. Note that all other sums and products are determined by the semiring axioms. This definition extends to the semiring $\T$ whose underlying set are the non-negative numbers $\R_{\geq0}$, whose multiplication is defined as usual and addition is defined as the maximum, i.e.\
 \[
  a+b \ = \ \max\{a,b\}.
 \]
 This semiring is called the \emph{tropical semifield}. 
 
 Note that the tropical semifield is typically described differently in the literature about tropical geometry. Namely, it takes the shape of the ``max-plus-algebra'' $\R\cup\{-\infty\}$ or the ``min-plus-algebra'' $\R\cup\{\infty\}$. However, all the three semirings, $\T$, $\R\cup\{-\infty\}$ and $\R\cup\{\infty\}$, are isomorphic. For our purposes, it is least confusing to use the non-negative real numbers with usual multiplication as a model for the tropical numbers.
\end{ex}

Given a semiring $R$, we can construct the ring $R_\Z=R\otimes_\N\Z$ of the formal differences of elements of $R$. In more detail, $R_\Z$ is defined as follows. As a set, $\R_\Z=R\times R/\sim$, where $(x,y)\sim (x',y')$ if and only if $x+y'+z=x'+y+z$ for some $z\in R$. We write $x-y$ for the equivalence class of $(x,y)$ in $R_\Z$. We can extend the addition and multiplication from $R$ to $R_\Z$ by the formulas
\[
 (x-y)+(z-t) = (x+z)-(y+t) \quad \text{and} \quad (x-y)\cdot(z-t) = (xz+yt)-(xt+yz).
\]

\begin{ex}
 We have $\N\otimes_\N\Z=\Z$, $\R_{\geq0}\otimes_\N\Z=\R$, $R[T_1,\dotsc,T_n]_\Z=R_\Z[T_1,\dotsc,T_n]$ and $\B\otimes_\N\Z=\T\otimes_\N\Z=\{0\}$.
\end{ex}

We can extend the notion of an ideal from rings to semirings: an \emph{ideal} of a semiring $R$ is a subset $I$ that contains $0$, $x+y$ and $tx$ for all $x,y\in I$ and $t\in R$. 

In contrast to ring theory, there is, however, no perfect unison between ideals and quotient semirings. This leads to the notion of a \emph{congruence}, which is an equivalence relation $\cR$ on a semiring $R$ such that the quotient set $R/\cR$ is a semiring, i.e.\ we can define addition and multiplication on equivalence classes by evaluating unambiguously on representatives. More explicitly, a congruence is an equivalence relation $\cR$ on $R$ such that $x\sim y$ and $z\sim t$ imply
\[
 x+z \ \sim \ y+t \qquad \text{and} \qquad xz \ \sim yt
\]
for all $x,y,z,t\in R$. For any set of relations $x_i\sim y_i$ between elements $x_i$ and $y_i$ of a semiring $R$, there is a smallest congruence $\cR$ containing these relations. Thus we can define
\[
 R\,/\,\gen{x_i\sim y_i} \ = \ R\,/\,\sim.
\]
Conversely, every morphism $f:R_1\to R_2$ of semirings has a \emph{congruence kernel}, which is the congruence $\cR$ on $R_1$ that is defined by $x\sim y$ if and only if $f(x)=f(y)$. In contrast, the \emph{kernel of $f$} is the ideal $\ker f=\{x\in R|f(x)=0\}$ of $R_1$.

Thus we gain two opposing constructions: we can associate with a congruence $\cR$ on $R$ the kernel of the quotient map $R\to R/\cR$, and we can associate with an ideal $I$ of $R$ the congruence $\cR=\gen{x\sim 0|x\in I}$. We write $R/I$ for $R/\cR$ in this case.

However, returning to our earlier remark, in general it is neither true that every congruence comes from an ideal nor that every ideal comes from a congruence. We call an ideal that occurs as the kernel of a morphism a \emph{$k$-ideal}. Note that if $R$ is a ring, then every ideal is a $k$-ideal.

\begin{ex}
 Every semiring can be written as a quotient of a polynomial semiring over $\N$, possibly in infinitely many indeterminants. For instance, we can define
 \[
  R \ = \ \N[T_1,\dotsc,T_4]\,/\,\gen{T_1T_4\sim T_2T_3+1}.
 \]
 Then $R_\Z=\Z[T_1,\dotsc,T_4]/(T_1T_4-T_2T_3-1)$ is the coordinate ring of $\SL(2)_\Z$.
\end{ex}

\subsubsection{Blueprints}
A \emph{blueprint} is a pair $B=(R,A)$ of a semiring $R$ and a multiplicative subset $A$ of $R$ that contains $0$ and $1$ and that generates $R$ as a semiring. 

The last condition in this definition is equivalent with the fact that $R$ is a quotient of the monoid semiring $\N[A]$ by a congruence $\cR$ on $\N[A]$. We also write $\bpquot A\cR$ for the blueprint $B=(R,A)$. Note that $R=\N[A]/\cR$ is determined by $A$ and $\cR$. Further, we write $B^\bullet=A$ and $B^+=R$. We write $x\= y$ if $(x,y)\in\cR$.

A \emph{blueprint morphism $f:B_1\to B_2$} is a semiring morphism $f^+:B_1^+\to B_2^+$ such that $f^+(B_1^\bullet)\subset B_2^\bullet$. Note that the restriction $f^\bullet:B_1^\bullet\to B_2^\bullet$ of $f^+$ is necessarily a monoid morphism and that $f^+$ is uniquely determined by $f^\bullet$.

\begin{ex}
 A monoid $A$ can be identified with the blueprint $(A,\N[A])$ and a semiring $R$ can be identified with the blueprint $(R,R)$. In this sense, blueprints are a simultaneous generalization of monoids and semirings. Note that if $B=(A,\N[A])$, then the base extension functor $-\otimes_\Fun\Z$ from section \ref{subsection: monoid schemes} can be recovered in terms of the formula
 \[
  A\otimes_\Fun\Z \ = \ B^+\otimes_\N\Z.
 \]
 
 There are also a series of novel constructions for blueprints. One of these are the cyclotomic field extensions $\F_{1^n}$ of $\Fun$, which are defined as follows for $n\geq 2$. Let $\Z[\zeta_n]$ be the ring of integers in the cyclotomic number field $\Q[\zeta_n]$ generated by a primitive $n$-th root $\zeta_n$ of unity and let $\mu_{n,0}=\{0\}\cup\{\zeta_n^i|i\in\Z\}$ be the submonoid generated by $\zeta_n$. Then $\F_{1^n}=(\mu_{n,0},\Z[\zeta_n])$ is a blueprint. It incorporates certain properties of the cyclotomic field $\Q(\zeta_n)$: we have $\F_{1^n}^+=\Z[\zeta_n]$ and the Galois group of $\Q(\zeta_n)/\Q$ equals the group of automorphisms of $\F_{1^n}$ that fix $\Fun$.
 
 Of particular importance is the case $n=2$: the ``quadratic'' extension $\Funsq=\bpgenquot{\{0,1,-1\}}{1+(-1)\=0}$ of $\Fun$ contains an additive inverse of $1$, which is important for several purposes. For instance, it is possible to define sheaf cohomology over $\Funsq$, cf.\ \cite{Flores-L-Szczesny17}. The blueprint $\Funsq$ also plays a role in our approach to algebraic groups over $\Fun$, see section \ref{subsubsection: the rank space}.
 
 Most interestingly for our purposes is that it is possible to define the blueprint
 \[
  B \ = \ \bpgenquot{\Fun[T_1,\dotsc,T_4]}{T_1T_4\= T_2T_3+1}.
 \]
 The coordinate ring of $\SL(2)_\Z$ can be recovered from this blueprints as $B^+_\Z$.
\end{ex}

\subsubsection{The spectrum}
Let $B$ be a blueprint. We will not enter all details in the definition of the spectrum $\Spec B$ of $B$, but concentrate on the description of the topological space associated with $\Spec B$. In fact, one can associate several topological spaces with $\Spec B$, cf.\ \cite{L15b} for more details. The relevant one for a theory of algebraic groups over $\Fun$ is based on the notion of a $k$-ideal, which is an ideal of the monoid $B^\bullet$ that spans a $k$-ideal in the semiring $B^+$.

More explicitly, a \emph{$k$-ideal} of a blueprint $B$ is a subset $I$ of $B^\bullet$ such that $0\in I$, $ab\in I$ for all $a\in I$ and $b\in B^\bullet$ and $c\in I$ whenever there are $a_1,\dotsc,a_n,b_,\dotsc,b_m\in B^\bullet$ such that $\sum a_i+c=\sum b_j$ in $R$. A $k$-ideal $\fp$ of $B$ is \emph{prime} if $S=B^\bullet-\fp$ is a multiplicative subset of $B^\bullet$.

We define $\Spec B$ as the set of all prime ideals $\fp$ of $B$ together with the topology generated by the open subsets
\[
 U_h \ = \ \bigl\{\, \fp\in\Spec B \, \bigl| \, h\notin \fp \, \bigl\}
\]
where $h$ ranges through all elements of $B^\bullet$. It carries a structure sheaf in blueprints, but its definition is somewhat more involved and we omit these details from our account.

\begin{ex}\label{ex: spectra of blueprints}
 If $B=(A,\N[A])$ is a monoid, then the last condition in the definition of a $k$-ideal is automatically satisfied. Thus $\Spec B$ equals the monoid spectrum $\Spec A$. For an arbitrary blueprint $B=(A,R)$, the spectrum $\Spec B$ is a subspace of $\Spec A$, together with the subspace topology.
 
 If $B=(R,R)$ is a semiring, then a $k$-ideal of $B$ is the same as a $k$-ideal of $R$. In particular if $R$ is a ring, then the spectrum $\Spec B$ coincides with the usual spectrum $\Spec R$ of the ring $R$. For an arbitrary blueprint $B=(A,R)$, we have a continuous map $\Spec R_\Z\to\Spec B$, which is surjective if $R=R_\Z$ is a ring. As a by-product, we obtain the notion of the spectrum of a semiring.
 
 As the next case, we inspect the spectrum of our leading example, $B=\bpgenquot{\Fun[T_1,\dotsc,T_4]}{T_1T_4\= T_2T_3+1}$. As explained above, it is a subspace of $\Spec\Fun[T_1,\dotsc,T_4]$, whose prime ideals are of the from $\fp_I=(T_i)_{i\in I}$ for any subset $I$ of $\{1,\dotsc,4\}$, cf.\ Example \ref{ex: monoid spectra}. In order to determine $\Spec B$, we have to verify, for which $I$, the set $\fp_I$ satisfies the additive axiom of a prime ideal of $B$. This can be tested on the generators of the defining congruence of $B$, which is $T_1T_4\= T_2T_3+1$.
 
 This relation implies that if both $T_1T_4$ and $T_2T_3$ are in a prime ideal $\fp$, then $1\in\fp$, which is not the case for any prime ideal. Thus we conclude that either $T_1T_4$ or $T_2T_3$ is not in $\fp$, which means that either $T_1$ and $T_4$ are not in $\fp$ or $T_2$ and $T_3$ are not in $\fp$. Thus $\Spec B$ consists of the prime ideals $\{0\}$, $(T_1),\dotsc,(T_4)$, $(T_1,T_4)$ and $(T_2,T_3)$. We illustrate $\Spec B$ in Figure \ref{figure: Spec of SL2 over Fun}.
  \begin{figure}[ht]
  \[
   \beginpgfgraphicnamed{fig7}
   \begin{tikzpicture}[inner sep=0,x=30pt,y=20pt]
    \draw (3,0) -- (0,2) -- (1,4) -- (2,2) -- (3,0) -- (4,2) -- (5,4) -- (6,2) -- (3,0);
    \filldraw (0,2) circle (3pt);
    \filldraw (2,2) circle (3pt);
    \filldraw (4,2) circle (3pt);
    \filldraw (6,2) circle (3pt);
    \filldraw (3,0) circle (3pt);
    \filldraw (1,4) circle (3pt);
    \filldraw (5,4) circle (3pt);
    \node at (3.8,0) {$\{0\}$};
    \node at (-0.5,2) {$(T_2)$};
    \node at (2.5,2) {$(T_3)$};
    \node at (3.5,2) {$(T_1)$};
    \node at (6.5,2) {$(T_4)$};
    \node at (0.2,4) {$(T_2,T_3)$};
    \node at (5.8,4) {$(T_1,T_4)$};
   \end{tikzpicture}
   \endpgfgraphicnamed
  \] 
  \caption{The spectrum of $\bpgenquot{\Fun[T_1,\dotsc,T_4]}{T_1T_4\sim T_2T_3+1}$}
  \label{figure: Spec of SL2 over Fun}
 \end{figure}
 
 We draw the reader's attention to the fact that the two closed points $(T_2,T_3)$ and $(T_1,T_4)$ stay in bijection to the elements of the Weyl group $W$ of $\SL(2)$, whose two elements are the subgroup $T$ of diagonal matrices in $\SL(2)$ and the set of  antidiagonal matrices. We shall investigate this fact in more depth in section \ref{subsection: algebraic groups over F1}.
\end{ex}

\subsubsection{Blue schemes}
We omit a rigorous definition of blue schemes, which would deviate into too heavy technicalities for the flavour of this overview paper. To give a taste, a blue scheme can be seen as a topological space together with a structure sheaf in blueprints that is locally isomorphic to the spectrum of a blueprint. This approach recovers monoid schemes and usual schemes as special cases; moreover, it provides a notion of semiring schemes.

We remark that there are other types of blue schemes, which are based on other types of prime ideals. These alternate approaches are relevant for different applications.

Meaningful variants are the following. The notion of a $k$-ideal, as considered above, yields blue schemes as introduced in \cite{L12a}. The notion of an ideal of the underlying monoid $B^\bullet$ of a blueprint $B$ yields blue schemes that were dubbed \emph{subcanonical} in \cite{L17}. Another variant is based on congruences. This has not worked out in the full generality of blueprints, but for certain subcategories by Berkovich (\cite{Berkovich11}) and Deitmar (\cite{Deitmar11a}).


\section{Growing up: achievements of \texorpdfstring{$\Fun$}{F1}-geometry}
\label{section: growing up}

Most of the initial goals of $\Fun$-geometry have been solved, with the exception of the most influential one, the Riemann hypothesis. 

The abc-conjecture has been claimed to be proven by Mochizuki in his monumental work \cite{Mochizuki12}. There are claims that some mathematicians have verified all details of the proof, but the acceptance by the community at large is still not clear at the time of writing. Although Mochizuki's proof follows a different line of thought, it contains ideas from $\Fun$-geometry, cf.\ \cite[Remark III.3.12.4 (iii)]{Mochizuki12} and \cite[Remark 5.10.2 (iii)]{Mochizuki15}.

$K$-theory has been developed for monoids and monoid schemes by Deitmar in \cite{Deitmar07} and by the author's collaboration \cite{Chu-L-Santhanam12} with Chu and Santhanam, respectively, and it has been shown that the $K$-theory of $\Fun$ coincides with the stable homotopy groups of the sphere.

Tits' dream of algebraic groups over $\Fun$ has been realized by the author in \cite{L12d}, based on the blueprint approach to $\Fun$. 

Besides settling these old scores, $\Fun$-geometry has found some further applications. Monoid schemes have been utilized by Corti\~nas, Haesemeyer, Walker and Weibel (\cite{CHWW15}) to connect the algebraic $K$-theory of toric varieties to cyclic homology in characteristic $p$.

The Giansiracusa brothers have used in \cite{Giansiracusa-Giansiracusa16} $\Fun$-schemes and semiring schemes to describe the tropicalization of a classical variety as a \emph{tropical scheme}. This approach towards tropical geometry promises to be a major breakthrough in the field. 

In the following, we explain some of these applications in more detail. In section \ref{subsection: steps towards the Riemann hypothesis}, we devote a few words to the progress that has been made towards the Riemann hypothesis. In section \ref{subsection: algebraic groups over F1}, we explain in a certain depth how blue schemes can be used to realize Tits’ dream of algebraic groups over $\Fun$. Though we try to make this accessible to the general reader, we have to assume a certain familiarity with group schemes when stating our main results. In section \ref{subsection: applications to tropical geometry}, we describe the relevance of tropical scheme theory for tropical geometry.


\subsection{Steps towards the Riemann hypothesis}\label{subsection: steps towards the Riemann hypothesis}

Several authors have considered compactifications $\overZ$ of $\Spec\Z$ and the arithmetic surface $\overZ\times_\Fun\overZ$; for instance, see \cite{Durov07}, \cite{Haran07}, \cite{L14b}, \cite{Takagi12b}. However, none of these ideas have been pursued further to the authors knowledge.

A somewhat different route, employing idempotent semirings, is taken by Connes and Consani who follow an ambitious programme around the Riemann hypothesis. Their research has already lead to a large number of publications; to name a few, cf.\ \cite{Connes-Consani15c}, \cite{Connes-Consani16c}, \cite{Connes-Consani17a}. We are not attempting an outline of this programme, but refer the interested reader to Connes' chapter in \cite{Nash-Rassias16} for such a summary.


\subsection{Algebraic groups over \texorpdfstring{$\Fun$}{F1}}\label{subsection: algebraic groups over F1}
In this section, we explain how we can make sense of Tits' dream of algebraic groups over $\Fun$ in the language of blueprints and blue schemes. Let $k$ be a ring and $G$ be a \emph{group scheme over $k$}, by which we mean a $k$-scheme $G$ together with a \emph{group law} $\mu:G\times G\to G$, a \emph{unit} $\epsilon:\Spec\Z\to G$ and an \emph{inversion} $\iota:G\to G$, which are $k$-linear morphisms that satisfy the usual axioms of a group, i.e.\ the diagrams
 \[
  \beginpgfgraphicnamed{fig8}
   \begin{tikzcd}[column sep=0.75cm]
    G\times G \times G\ar{r}{\id\times\mu} \ar{d}{\mu\times\id} & G\times G \ar{d}{\mu} &                                          & G\times G \ar{d}{\mu} & G \ar{r}{\Delta} \ar{d} & G\times G \ar{r}{\id\times\iota} & G\times G \ar{d}{\mu}  \\
    G\times G \ar{r}{\mu}                                       & G                     & G \ar{r}{\id} \ar{ru}{\id\times\epsilon} & G                     & \Spec k \ar{rr}{\epsilon} &                                & G                
   \end{tikzcd}
  \endpgfgraphicnamed
 \] 
commute where $\Delta:G\to G\times G$ is the diagonal and the $G\to\Spec k$ is the unique morphism to $\Spec k$. As it is the case for usual groups, the unit $\epsilon$ and the inversion $\iota$ are uniquely determined by $G$ and $\mu$.


\begin{ex}\label{ex: group law of SL_2Z}
 As our main example that guides us through the concepts of this section, we will inspect the special linear group $\SL(2)$. As a scheme, $\SL(2)_\Z$ is defined as the spectrum of $\Z[T_1,\dotsc,T_4]/(T_1T_4-T_2T_3-1)$. We have already seen in Example \ref{ex: spectra of blueprints} that $SL(2)_\Z$ descends to $\Fun$ as a blue scheme.
 
 The group law $\mu_\Z:\SL(2)_\Z\times\SL(2)_\Z\to \SL(2)_\Z$ is determined by the matrix multiplication
 \[
  \begin{pmatrix} T_1 & T_2 \\ T_3 & T_4 \end{pmatrix} \ = \ 
  \begin{pmatrix} T_1' & T_2' \\ T_3' & T_4' \end{pmatrix} \cdot
  \begin{pmatrix} T_1'' & T_2'' \\ T_3'' & T_4'' \end{pmatrix} \ = \ 
  \begin{pmatrix} T_1'T_1''+T_2'T_3'' & T_1'T_2''+T_2'T_4'' \\ T_3'T_1''+T_4'T_3'' & T_3'T_2''+T_4'T_4'' \end{pmatrix}
 \]
 of $2\times2$-matrices with determinant $1$.
\end{ex}

To put Tits' idea of algebraic groups over $\Fun$ into modern language, we ask for the following: given a group scheme $G$ (over $\Z$) together with a group law $\mu:G\times G\to G$, is there an $\Fun$-scheme $G_\Fun$ together with a group law $\mu_\Fun:G_\Fun\times G_\Fun\to G_\Fun$ whose base extensions to usual schemes are $G$ and $\mu$, respectively?

While some approaches to $\Fun$-geometry contain $\Fun$-models $G_\Fun$ for large classes of group schemes $G$ it is a problem to descend the group law to a morphism of $\Fun$-schemes. Heuristically speaking, the group law of most group schemes involves the addition of coordinates, as it is the case for $\SL(2)_\Z$, cf.\ Example \ref{ex: group law of SL_2Z}. Therefore it is not possible to consider such a morphism in a \emph{non-additive} geometry. 

There is a second, more subtle, problem concerning the formula $G(\Fun)=W$ where $G(\Fun)=\Hom(\Spec\Fun,G)$ and where $W$ is the Weyl group of $G$. Namely, if this equality was an identity of groups, then it would imply that the Weyl group $W$ of $G$ could be embedded as a subgroup of $G(\Z)$ such that each element of $W$ lies in the coset representing it. But this leads already in the simplest case of $G=\SL(2)$ to a contradiction: its Weyl group is $\{\pm1\}$, but the coset of $-1$ is the anti-diagonal in $\SL(2,\Z)$, which does not contain an element of order $2$. For more details on this, cf.\ \cite{L16}.

Our conclusion is the following: we are required to alter the notion of a morphism to make it possible to descend group laws and we cannot take the formula $G(\Fun)=W$ literally, but have to find an appropriate meaning for it. This has been done in \cite{L12b}, using torified schemes and $\Fun$-schemes after Connes and Consani, and in \cite{L12d}, using blue schemes. We describe the latter approach in the upcoming paragraphs.

\subsubsection{Algebraic tori}
 
Because of their central role for what is to come we begin with a description of (split) algebraic tori and their $\Fun$-models. Algebraic tori are the only connected group schemes that fit easily into any concept of $\Fun$-geometry. We examine the case of a torus of rank $1$ from the perspective of blue schemes; the case of higher rank can be deduced easily from the following description.

Let $B$ be a blueprint. We define the \emph{multiplicative group scheme $\G_{m,B}$ over $B$} as $\Spec B[T^{\pm1}]$. The multiplication $\mu:\G_{m,B}\times \G_{m,B}\to\G_{m,B}$ is given by the morphism
\[
 \begin{array}{cccc}
  \mu^\#: & B[T^{\pm1}] & \longrightarrow & B[T^{\pm1}] \otimes_B B[T^{\pm1}], \\
              & aT       & \longmapsto     & aT\otimes T
 \end{array}
\]
the unit $\epsilon:\Spec B\to \G_{m,B}$ is given by the morphism
\[
 \begin{array}{cccc}
  \epsilon^\#: & B[T^{\pm1}] & \longrightarrow & B \\
              & aT       & \longmapsto     & a
 \end{array}
\]
and the inversion $\iota:\G_{m,B}\to \G_{m,B}$ is given by the morphism
\[
 \begin{array}{cccc}
  \iota^\#: & B[T^{\pm1}] & \longrightarrow & B[T^{\pm1}]. \\
              & aT       & \longmapsto     & aT^{-1}
 \end{array}
\]
If $B$ is $\Fun$ or $\Funsq=\bpgenquot{\{0,\pm1\}}{1+(-1)\=0}$, then the base extension of $\G_{m,B}$ to $\Z$ is the multiplicative group scheme $\G_{m,\Z}$ over $\Z$ as a group scheme.\footnote{Please note that we face a clash of notation at this point: while we denote by $\G_{m,\Z}$ the spectrum of the \emph{polynomial ring} $\Z[T]^+$, the very same notation is also used for the spectrum of the free blueprint $(\{aT^i\},\Z[T]^+)$ in the definition of $\G_{m,B}$ in the case $B=\Z$. However, for the sake of a more intuitive notation, we do not dissolve this contradiction, but refer the reader to \cite{L12d} and \cite{L16} for a more sophisticated treatment.}

\subsubsection{The rank space}
\label{subsubsection: the rank space}

A blueprint $B$ is \emph{cancellative} if $\sum a_i +c\=\sum b_j+c$ implies $\sum a_i\=\sum b_j$. Note that a blueprint $B$ is cancellative if and only if the canonical morphism $B\to B^+_\Z$ is injective. A blue scheme is cancellative if its structure sheaf takes values in the class of cancellative blueprints. 

Let $X$ be a blue scheme and $x\in X$ a point. As in usual scheme theory, every closed subset of $X$ comes with a natural structure of a \emph{(reduced) closed subscheme of $X$}, cf.\ \cite[Section 1.4]{L12d}. The \emph{rank $\rk\, x$ of $x$} is the dimension of the $\Q$-scheme $\barx^+_\Q$ where $\barx$ denotes the closure of $x$ in $X$ together with its natural structure as a closed subscheme. Define
\[
 r \quad = \quad \min\, \{\ \rk\, x \ | \ x\in X \ \}.
\]
For the sake of simplicity, we will make the following general hypothesis on $X$. 
\begin{enumerate}[label={(H)}]
 \item\label{hypothesisH} The blue scheme $X$ is connected and cancellative. For all $x\in X$ with $\rk\, x=r$, the closed subscheme $\barx$ of $X$ is isomorphic to either $\G_{m,\Fun}^{r}$ or $\G_{m,\Funsq}^{r}$.
\end{enumerate}
This hypothesis allows us to surpass certain technical aspects in the definition of the rank space. 

Assume that $X$ satisfies (H). Then the number $r$ is denoted by $\rk\,X$ and is called the \emph{rank of $X$}. The \emph{rank space of $X$} is the blue scheme
\[
 X^\rk\quad = \quad \coprod_{\rk\, x=r}\ \barx,
\]
and it comes together with a closed immersion $\rho_X:X^\rk\to X$. 

Note that Hypothesis \ref{hypothesisH} implies that $X^\rk$ is the disjoint union of copies of tori $\G_{m,\Fun}^r$ and $\G_{m,\Funsq}^r$. Since the underlying set of both $\G_{m,\Fun}^r$ and $\G_{m,\Funsq}^r$ is the one-point set, the underlying set of $X^\rk$ is $\cW(X)=\{x\in X|\rk\,x=r\}$. 

For a blue scheme $Y$, we denote by $\beta_Y:Y^+_\Z\to Y$ the base extension morphism. We have $X^{\rk,+}_\Z\simeq\coprod_{\rk\,x=r} \G_{m,\Z}^r$ and we obtain a commutative diagram
\[
 \xymatrix@C=4pc{   X^{\rk,+}_\Z \ar[r]^{\rho^+_{X,\Z}}\ar[d]^{\beta_{X^\rk}}   &  X^{+}_\Z \ar[d]^{\beta_{X}} \\
                    X^\rk \ar[r]^{\rho_X}                                    &  X\ .\hspace{-5pt} } 
\]

\begin{ex}\label{ex: rank space of SL_2}
 We determine the rank space of our leading example $\SL(2)$ and its $\Fun$-model $\SL(2)_\Fun$. It is not hard to verify that $\SL(2)_\Fun$ is cancellative. As already explained in Example \ref{ex: spectra of blueprints}, the points of $\SL(2)_\Fun$ are of the form $(T_i)_{i\in I}$ where $I$ is a subset of $\{1,2,3,4\}$ that does not contain elements of both $\{2,3\}$ and $\{1,4\}$. The subscheme $\big(\overline{(T_i)_{i\in I}}\big)^+_\Q$ of $\SL(2)_\Q$ represents the set of $2\times 2$-matrices ${\scriptstyle\begin{psmallmatrix}T_1& T_2\\ T_3& T_4\end{psmallmatrix}}$ with determinant $1$ for which $T_i=0$ for all $i\in I$. For instance, $\overline{(T_2,T_3)}^+_\Q$ is the diagonal torus in $\SL(2)_\Q$, which is one dimensional, $\overline{(T_1,T_4)}^+_\Q$ is the antidiagonal torus, also one dimensional, $\overline{(T_2)}^+_\Q$ is the Borel subgroup of upper triangular matrices, which is of dimension $2$, and $\overline{\{0\}}^+_\Q$ equals $\SL(2)_\Q$. In Figure \ref{figure: ranks of points of SL_2}, we illustrate all points of $\SL(2)_\Fun$, together with their ranks.
 
 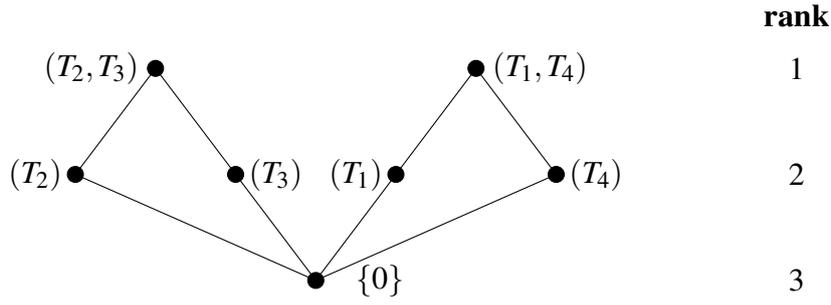
\begin{figure}
  \[
   \beginpgfgraphicnamed{fig9}
    \begin{tikzpicture}[inner sep=0,x=30pt,y=20pt]
     \draw (3,0) -- (0,2) -- (1,4) -- (2,2) -- (3,0) -- (4,2) -- (5,4) -- (6,2) -- (3,0);
     \filldraw (0,2) circle (3pt);
     \filldraw (2,2) circle (3pt);
     \filldraw (4,2) circle (3pt);
     \filldraw (6,2) circle (3pt);
     \filldraw (3,0) circle (3pt);
     \filldraw (1,4) circle (3pt);
     \filldraw (5,4) circle (3pt);
     \node at (3.8,0) {$\{0\}$};
     \node at (-0.5,2) {$(T_2)$};
     \node at (2.5,2) {$(T_3)$};
     \node at (3.5,2) {$(T_1)$};
     \node at (6.5,2) {$(T_4)$};
     \node at (0.2,4) {$(T_2,T_3)$};
     \node at (5.8,4) {$(T_1,T_4)$};
     \node at (9,5) {\textbf{rank}};
     \node at (9,4) {$1$};
     \node at (9,2) {$2$};
     \node at (9,0) {$3$};
    \end{tikzpicture}
   \endpgfgraphicnamed
  \]
  \caption{The ranks of the points of $\SL(2)_\Fun$}
  \label{figure: ranks of points of SL_2}
 \end{figure}

 We conclude that $(T_2,T_3)$ and $(T_1,T_4)$ are the points of minimal rank. Since
 \begin{align*}
  \overline{(T_2,T_3)} \ &= \ \Spec\big(\bpgenquot{\Fun[T_1,T_2,T_3,T_4]}{T_1T_4\=T_2T_3+1,T_2=T_3=0}\big) \\ &= \ \Spec\big(\bpgenquot{\Fun[T_1,T_4]}{T_1T_4\=1}\big) \ &&\simeq \ \G_{m,\Fun}
 \end{align*}
 and 
 \begin{align*}
  \overline{(T_1,T_4)} \ &= \ \Spec\big(\bpgenquot{\Fun[T_1,T_2,T_3,T_4]}{T_1T_4\=T_2T_3+1,T_1=T_4=0}\big) \\ &= \ \Spec\big(\bpgenquot{\Fun[T_2,T_3]}{T_2T_3+1\=0}\big) \ && \simeq \ \G_{m,\Funsq},
 \end{align*}
 we conclude that hypothesis \ref{hypothesisH} is satisfied by $\SL(2)_\Fun$ and its rank space is
 \[
  \SL(2)_\Fun^\rk \ = \ \overline{(T_2,T_3)} \ \amalg \ \overline{(T_1,T_4)} \ \simeq \ \G_{m,\Fun} \ \amalg \ \G_{m,\Funsq}.
 \]
\end{ex}

\subsubsection{The Tits category and the Weyl extension}
\label{subsubsection: The Tits category and the Weyl extension}

A \emph{Tits morphism $\varphi:X\to Y$} between two blue schemes $X$ and $Y$ is a pair $\varphi=(\varphi^\rk,\varphi^+)$ of a morphism $\varphi^\rk: X^\rk\to Y^\rk$ and a morphism $\varphi^+: X^+\to Y^+$ such that the diagram
\[
 \begin{tikzcd}[column sep=1.5cm]
  X^{\rk,+}_\Z  \ar{r}{\varphi^{\rk,+}_\Z} \ar{d}[left]{\rho_{X,\Z}^+}      & Y^{\rk,+}_\Z  \ar{d}{\rho_{Y,\Z}^+} \\ 
           X^+_\Z  \ar{r}{\varphi^+_\Z}                                      & Y^+_\Z  
 \end{tikzcd}
\]
commutes where $\varphi^+_\Z$ and $\varphi^{\rk,+}_\Z$ are the base extension of $\varphi^+$ and $\varphi^\rk$, respectively, to usual schemes. We denote the category of blue schemes together with Tits morphisms by $\Sch_\cT$ and call it the \emph{Tits category}.

The Tits category comes together with two important functors. The \emph{Weyl extension $\cW:\Sch_\cT\to\Sets$} sends a blue scheme $X$ to the underlying set $\cW(X)=\{x\in X|\rk\,x=r\}$ of $X^\rk$ and a Tits morphism $\varphi:X\to Y$ to the underlying map $\cW(\varphi):\cW(X)\to \cW(Y)$ of the morphism $\varphi^\rk:X^\rk\to Y^\rk$. 

The base extension $(-)^+:\Sch_\cT\to\Sch^+$ sends a blue scheme $X$ to its universal semiring scheme $X^+$ and a Tits morphism $\varphi:X\to Y$ to $\varphi^+:X^+\to Y^+$. We obtain the following diagram of ``base extension functors''
\[
 \begin{tikzcd}[column sep=1.5cm]
  \Sets & & \Sch_\Z & & \Sch_R \\
        & &         & \Sch_\N \ar[-]{ul}{(-)_\Z^+}\ar[-]{ur}[swap]{(-)_R^+} \\
                 & & \TSch \ar[-]{uull}{\cW}\ar[-]{ur}[swap]{(-)^+} 
 \end{tikzcd}
\]
from the Tits category $\TSch$ to the category $\Sets$ of sets, to the category $\Sch^+_\Z$ of usual schemes and to the category $\Sch^+_R$ of semiring schemes over any semiring $R$.

\begin{thm}[{\cite[Thm.\ 3.8]{L12d}}]
 All functors appearing in the above diagram commute with finite products. Consequently, all functors send (semi)group objects to (semi)group objects.
\end{thm}

\begin{ex}\label{ex: descending the group law of SL_2 to F1}
 The group law $\mu_\Z$ of $\SL(2)_\Z$ descends uniquely to a Tits morphism 
 \[
  \mu_\Fun=(\mu^+,\mu^\rk):\SL(2)_\Fun\times\SL(2)_\Fun\to\SL(2)_\Fun
 \]
 in $\Sch_\cT$. To wit, $\mu_\Z$ descends uniquely to a morphism $\mu^+:\SL(2)_\N\times\SL(2)_\N\to\SL(2)_\N$ where 
 \[
  \SL(2)_\N=\SL(2)_\Fun^+=\Spec\big(\N[T_1,T_2,T_3,T_4]/(T_1T_4\sim T_2T_3+1)\big)
 \]
 since $\mu_\Z$ is defined by means of sums and products, without employing subtraction. The morphism $\mu^\rk$ is determined by the commutativity of the diagram 
 \[
  \begin{tikzcd}[column sep=2.0cm]
   \big(SL(2)_\Fun^\rk \times SL(2)_\Fun^\rk\big)^+_\Z \ar{r}{\mu^{\rk,+}_\Z} \ar{d}[left]{}      & \big(SL(2)_\Fun^\rk\big)^+_\Z \ar{d}{} \\ 
            \big(\SL(2)_\N\times\SL(2)_\N\big)_\Z  \ar{r}{\mu^+_\Z}                                      & \big(\SL(2)_\N\big)_\Z.  
  \end{tikzcd}
 \]
 The Weyl extension of $\SL(2)_\Fun$ consists of two points $x_{23}$ and $x_{14}$, which are the respective unique points of the components $\overline{(T_2,T_3)}$ and $\overline{(T_1,T_4)}$ of the rank space $\SL(2)_\Fun^\rk$. The Weyl extension of $\mu_\Fun$ endows $\cW\big(\SL(2)_\Fun\big)$ with the structure of a group with neutral element $x_{23}$.
\end{ex}

\subsubsection{Tits-Weyl models}
\label{subsubsection: Tits-Weyl models}

A \emph{Tits monoid} is a (not necessarily commutative) monoid in $\Sch_\cT$, i.e.\ a blue scheme $G$ together with an associative multiplication $\mu:G\times G\to G$ in $\Sch_\cT$ that has an identity $\epsilon:\Spec\Fun\to G$. We often understand the multiplication $\mu$ implicitly and refer to a Tits monoid simply by $G$. The Weyl extension $\cW(G)$ of a Tits monoid $G$ is a unital associative semigroup. The base extension $G^+$ is a (not necessarily commutative) monoid in $\Sch_\N$.

Given a Tits monoid $G$ satisfying (H) with multiplication $\mu$ and identity $\epsilon$, the image of $\epsilon:\Spec\Fun\to G$ consists of a closed point $e$ of $X$. The closed reduced subscheme $\fe=\{e\}$ of $G$ is called the \emph{Weyl kernel of $G$}. 

\begin{lemma}[{\cite[Lemma 3.11]{L12d}}]
 The multiplication $\mu$ restricts to $\fe$, and with this, $\fe$ is isomorphic to the torus $\G_{m,\Fun}^r$ as a group scheme where $r=\rk X$. 
\end{lemma}

This means that $\fe^+_\Z$ is a split torus $T\simeq\G_{m,\Z}^r$ of $\cG=G^+_\Z$, which we call the \emph{canonical torus of $\cG$ (w.r.t.\ $G$)}. 

If $\cG$ is an affine smooth group scheme of finite type, then we obtain a canonical morphism 
\[
 \Psi_\fe : \ G^{\rk,+}_\Z/\fe^+_\Z \quad \longrightarrow \quad W(T)
\]
of schemes where $W(T)=\Norm_\cG(T)/\Cent_\cG(T)$ is the \emph{Weyl group of $\cG$ w.r.t.\ $T$} (cf.\ \cite[XIX.6]{SGA3III}). We say that $G$ is a \emph{Tits-Weyl model of $\cG$} if $T$ is a maximal torus of $\cG$ (cf.\ \cite[XII.1.3]{SGA3II}) and $\Psi_\fe$ is an isomorphism. 

\begin{ex}\label{Tits-Weyl model of SL_2}
 The blue scheme $\SL(2)_\Fun$ together with the Tits morphism $\mu_\Fun$ is a Tits-Weyl model of $\SL(2)_\Z$, which can be reasoned as follows. It can be easily seen that $\mu_\Fun$ is an associative multiplication for $\SL(2)_\Fun$ in the Tits category, which has a unit $\epsilon:\Spec\Fun\to \SL(2)_\Fun$, which is the pair of the morphism $\epsilon^+:\Spec\N\to \SL(2)_\N$, given by $\scriptstyle\begin{psmallmatrix} T_1 & T_2 \\ T_3 & T_4 \end{psmallmatrix} \mapsto \begin{psmallmatrix} 1 & 0 \\ 0 & 1 \end{psmallmatrix}$, and the morphism $\epsilon^\rk:\Spec \Fun\to \overline{(T_2,T_3)}\simeq\Spec\big(\Fun[T_1^{\pm1}]\big)$, given by $T_1\mapsto 1$. 
 
 The Weyl kernel $\fe$ of the Tits monoid $\SL(2)_\Fun$ is therefore $\overline{(T_2,T_3)}\simeq\Spec\big(\Fun[T_1^{\pm1}]\big)$, and we have
 \[
  \fe^+_\Z \ = \ \overline{(T_2,T_3)}^+_\Z \ = \ \Spec\big(\Z[T_1,T_4]/(T_1T_4-1)\big) \ = \ T \ = \ \Cent_{\SL(2)_\Z}(T)
 \]
 where the canonical torus $T$ is the diagonal torus, which is a maximal torus of $\SL(2)_\Z$. Since the normalizer of $T$ in $\SL(2)_\Z$ is the union of the diagonal torus $T$ with the antidiagonal torus, we obtain an isomorphism
 \[
  \Psi_\fe: \ \big(\SL(2)_\Fun\big)^{\rk,+}_\Z/\fe^+_\Z \quad \stackrel\sim\longrightarrow \quad \Norm_{\SL(2)_\Z}(T)/\Cent_{\SL(2)_\Z}(T) \ = \ W(T).
 \]
 This shows that $\SL(2)_\Fun$ together with $\mu_\Fun$ is a Tits-Weyl model of $\SL(2)_\Z$.

\end{ex}

We review some definitions, before we formulate the properties of Tits-Weyl models in Theorem \ref{thm: properties of tits-weyl groups} below. The \emph{ordinary Weyl group of $\cG$} is the underlying group $W$ of $W(T)$. The \emph{reductive rank of $\cG$} is the rank of a maximal torus of $\cG$. For a split reductive group scheme, we denote the \emph{extended Weyl group} or \emph{Tits group} $\Norm_\cG(T)(\Z)$ by $\widetilde W$ (cf.\ \cite{Tits66} or \cite[Section 3.3]{L12d}).

For a blueprint $B$, the set $G^\sT(B)$ of Tits morphisms from $\Spec B$ to $G$ inherits the structure of an associative unital semigroup. In case $G$ has several connected components, we define the \emph{rank of $G$} as the rank of the connected component of $G$ that contains the image of the unit $\epsilon:\Spec\Fun\to G$. 

\begin{thm}[{\cite[Thm.\ 3.14]{L12d}}]\label{thm: properties of tits-weyl groups}
 Let $\cG$ be an affine smooth group scheme of finite type. If $\cG$ has a Tits-Weyl model $G$, then the following properties hold true.
 \begin{enumerate}
  \item\label{part1} The Weyl group $\cW(G)$ is canonically isomorphic to the ordinary Weyl group $W$ of $\cG$.
  \item\label{part2} The rank of $G$ is equal to the reductive rank of $\cG$.
  \item\label{part3} The semigroup $\Hom_\cT(\Spec\Fun,G)$ is canonically a subgroup of $\cW(G)$.
  \item\label{part4} If $\cG$ is a split reductive group scheme, then $\Hom_\cT(\Spec\Funsq,G)$ is canonically isomorphic to the extended Weyl group $\widetilde W$ of $\cG$.
 \end{enumerate}
\end{thm}

The following theorem is proven in \cite{L12d} for a large class of split reductive group schemes $\cG$ and their Levi- and parabolic subgroups. An additional idea of Reineke extended this result to all split reductive group schemes; cf.\ \cite{L16}.

\begin{thm}\label{thm: chevalley groups over f1}
 \begin{enumerate}
  \item Every split reductive group scheme $\cG$ has a Tits-Weyl model $G$.
  \item Let $T$ be the canonical torus of $\cG$ and $\cM$ a Levi subgroup of $\cG$ containing $T$. Then $\cM$ has a Tits-Weyl model $M$ that comes together with a locally closed embedding $M\to G$ of Tits-monoids that is a Tits morphism.
  \item Let $\cP$ a parabolic subgroup of $\cG$ containing $T$. Then $\cP$ has a Tits-Weyl model $P$ that comes together with a locally closed embedding $P\to G$ of Tits-monoids that is a Tits morphism.
  \item Let $\cU$ be the unipotent radical of a parabolic subgroup $\cP$ of $\GL_{n,\Z}$ that contains the diagonal torus $T$. Then $\cU$, $\cP$ and $\GL_{n,\Z}$ have respective Tits-Weyl models $U$, $P$ and $\GL_{n,\Fun}$, together with locally closed embeddings $U\to P\to \GL_{n,\Fun}$ of Tits-monoids that are Tits morphisms and such that $T$ is the canonical torus of $\cP$ and $\GL_{n,\Z}$.
 \end{enumerate}
\end{thm}

\subsubsection{Tits' dream}
\label{subsubsection: Tits dream}

With the formalism developed in the previous sections, it is possible to make Tits' idea precise: the combinatorial counterparts of geometries over $\F_q$ can be seen as geometries over $\Fun$. This is explained in detail in \cite[section 6]{L16}, using the language of buildings.

We explain how this works in the previously considered example $\GL(3)$, cf.\ Example \ref{ex: Tits geometry of GL_3}. Similar to the example of $\SL(2)$, the group scheme $\GL(3)_\Z$ has a Tits-Weyl model $\GL(3)_\Fun$, which is a monoid in the Tits category. The Weyl extension of $\GL(3)_\Fun$ is the symmetric group $S_3$ on three elements, which equals the Weyl group of $\GL(3)$. 

The standard action of $\GL(3)_\Z$ on the Grassmannians $\Gr(k,3)_\Z$ for $k=1,2$ descends to monoid actions of $\GL(3)_\Fun$ on the standard $\Fun$-models $\Gr(k,3)_\Fun$ of the Grassmannians, considered as objects of the Tits category $\Sch_\cT$. The Weyl extension of these actions corresponds to the action of $S_3$ on the set $\Sigma(1,3)=\{1,2,3\}$ for $k=1$ and to the action of $S_3$ on the set $\Sigma(2,3)=\big\{\{1,2\},\{1,3\},\{2,3\}\big\}$ for $k=2$.

Replacing $\Hom(\Spec\Fun, X)$ by $\cW(X)$ makes precise the heuristics
\[
 \cW(\GL(3)_\Fun) \ = \ S_3 \qquad \text{and} \qquad \cW(\Gr(k,3)) \ = \ \Sigma(k,3).
\]
This example is generalized to all other classical groups in \cite[section 6]{L16}.


\subsection{Applications to tropical geometry}
\label{subsection: applications to tropical geometry}

An exciting new application of $\Fun$-geometry lies in the field of tropical geometry. Such an application was already foreseen in various works on $\Fun$-geometry, cf.\ \cite{Durov07}, \cite{Toen-Vaquie09} and \cite{L12a}, but it took until 2013 till Jeffrey and Noah Giansiracusa (\cite{Giansiracusa-Giansiracusa16}) used $\Fun$-schemes and semiring schemes to underlay tropical varieties with the structure of a semiring scheme---at least in the case of the tropicalization of a classical scheme.

We will not dwell into this theory, but restrict ourselves to an explanation of the relevance of a scheme theoretic approach for tropical geometry and a description of the first results in this direction. For more details, we refer to Giansiracusas' papers \cite{Giansiracusa-Giansiracusa14} and \cite{Giansiracusa-Giansiracusa16}, to Maclagan and Rinc\'on's papers \cite{Maclagan-Rincon14} and \cite{Maclagan-Rincon16} as well as to the author's paper \cite{L15b}.

\subsubsection{Some history}

Tropical geometry was born with Mikhalkin's calculation (\cite{Mikhalkin00}) of Gromov-Witten invariants around twenty years ago. He was able to convert the classical problem of counting the number of nodal algebraic curves passing through a given number of points to a counting problem of tropical curves, and to solve the latter problem by means of elementary combinatorics.

Since then tropical geometry has developed rapidly. While during the first decade of tropical geometry, researchers concentrated mainly on applications in the vein of Mikhalkin's results, tropical geometry entered a second era with the fundamental works of Kajiwara (\cite{Kajiwara08}) and Payne (\cite{Payne09}) around 10 years ago. They provided an elegant framework for tropical geometry and tied it to nonarchimedean analytic geometry in the sense of Berkovich. Soon after these insights, tropical geometry found applications in Brill-Noether theory, moduli spaces, skeleta of Berkovich spaces and rational points. 

\subsubsection{Tropical varieties}
In order to explain the relevance of tropical scheme theory, we have to explain what a tropical variety is. To this end, it is easiest to work with the max-plus algebras $\R\cup\{-\infty\}$ as the tropical numbers $\T$ where addition is the maximum, $a+^\trop b=\max\{a,b\}$, and its multiplication the usual addition, $a\ast^\trop b=a+b$, with the obvious extensions of these operations to $-\infty$; cf.\ Example \ref{ex: semirings} for alternative descriptions of the tropical numbers. In this form, the multiplicative units of $\T$ are $\T^\times=\R$.

In its simplest incarnation, a \emph{tropical variety} is a closed subset of $(\T^\times)^n=\R^n$ together with a subdivision into finitely many rational polyhedra and together with a weight function on the top dimensional polyhedra that satisfies a certain balancing condition. These definitions somewhat difficult to explain for higher dimensions---and we choose to omit them---, but in the case of dimension $1$, they boil down to the following.

A \emph{tropical curve in $\R^n$} is a finite graph with possibly unbounded edges, embedded in $\R^n$, together with a weight function in $\Z_{>0}$ on the edges such that all edges have rational slope and such that for every vertex $p$ the following balancing condition is satisfied: given an edge $e$ adjacent to $p$, let $v_e$ be primitive vector at $p$ pointing in the direction of $e$, i.e.\ $v_e$ is the shortest vector in $\Z^2-\{0\}$ such that $v_e$ is contained in the ray generated by $\{q-p|q\in e\}\subset\R^2$; then we have
\[
 \sum_{e\text{ adjacent to }p} v_e \ = \ 0.
\]

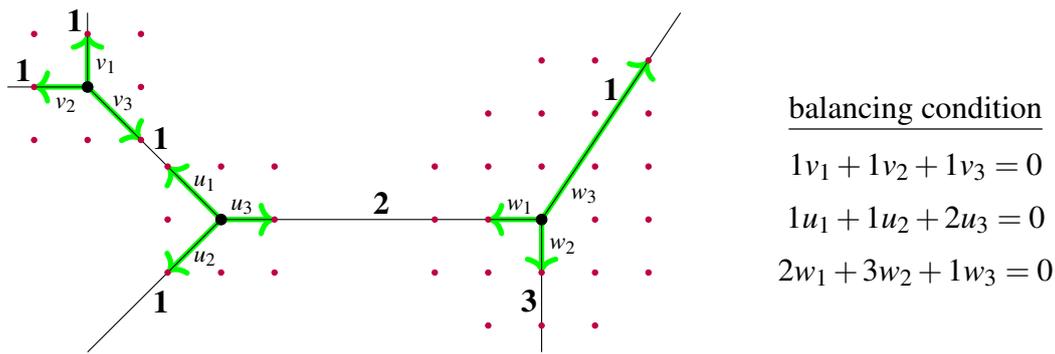
\begin{figure}[ht]
 \[
 \beginpgfgraphicnamed{fig11}
  \begin{tikzpicture}[inner sep=0,x=20pt,y=20pt]
   \draw[->,green,line width=2pt] (0.5,4.5)-- node[pos=0.4,right=2pt,black] {$\scriptstyle v_1$} (0.5,5.5);
   \draw[->,green,line width=2pt] (0.5,4.5)-- node[pos=0.4,below=3pt,black] {$\scriptstyle v_2$} (-0.5,4.5);
   \draw[->,green,line width=2pt] (0.5,4.5)-- node[pos=0.3,right=2.5pt,black] {$\scriptstyle v_3$} (1.5,3.5);
   \draw[->,green,line width=2pt] (3,2)-- node[pos=0.3,above=4pt,black] {$\scriptstyle u_1$} (2,3);
   \draw[->,green,line width=2pt] (3,2)-- node[pos=0.3,below=5pt,black] {$\scriptstyle u_2$} (2,1);
   \draw[->,green,line width=2pt] (3,2)-- node[pos=0.4,above=1pt,black] {$\scriptstyle u_3$} (4,2);
   \draw[->,green,line width=2pt] (9,2)-- node[pos=0.4,above=1pt,black] {$\scriptstyle w_1$} (8,2);
   \draw[->,green,line width=2pt] (9,2)-- node[pos=0.5,right=2pt,black] {$\scriptstyle w_2$} (9,1);
   \draw[->,green,line width=2pt] (9,2)-- node[pos=0.15,right=4pt,black] {$\scriptstyle w_3$} (11,5);
   \draw (0.5,4.5)-- node[pos=0.8,above=2pt,black] {$\bf 1$} (-1,4.5);
   \draw (0.5,4.5)-- node[pos=0.9,left=1.5pt,black] {$\bf 1$} (0.5,5.9);
   \draw (0.5,4.5)-- node[pos=0.55,above=4pt,black] {$\bf 1$} (3,2);
   \draw (0.5,-0.5)-- node[pos=0.55,below=4pt,black] {$\bf 1$} (3,2);
   \draw (3,2)-- node[pos=0.5,above=2pt,black] {$\bf 2$} (9,2);
   \draw (9,2)-- node[pos=0.5,above=6pt,black] {$\bf 1$} (11.6,5.9);
   \draw (9,-0.5)-- node[pos=0.39,left=1.5pt,black] {$\bf 3$} (9,2);
   \filldraw (0.5,4.5) circle (2pt);
   \filldraw[purple] (-0.5,3.5) circle (1pt);
   \filldraw[purple] (-0.5,4.5) circle (1pt);
   \filldraw[purple] (-0.5,5.5) circle (1pt);
   \filldraw[purple] (0.5,3.5) circle (1pt);
   \filldraw[purple] (0.5,5.5) circle (1pt);
   \filldraw[purple] (1.5,3.5) circle (1pt);
   \filldraw[purple] (1.5,4.5) circle (1pt);
   \filldraw[purple] (1.5,5.5) circle (1pt);
   \filldraw (3,2) circle (2pt);
   \filldraw[purple] (2,1) circle (1pt);
   \filldraw[purple] (2,2) circle (1pt);
   \filldraw[purple] (2,3) circle (1pt);
   \filldraw[purple] (3,1) circle (1pt);
   \filldraw[purple] (3,3) circle (1pt);
   \filldraw[purple] (4,1) circle (1pt);
   \filldraw[purple] (4,2) circle (1pt);
   \filldraw[purple] (4,3) circle (1pt);
   \filldraw (9,2) circle (2pt);
   \filldraw[purple] (7,1) circle (1pt);
   \filldraw[purple] (7,2) circle (1pt);
   \filldraw[purple] (7,3) circle (1pt);
   \filldraw[purple] (8,0) circle (1pt);
   \filldraw[purple] (8,1) circle (1pt);
   \filldraw[purple] (8,2) circle (1pt);
   \filldraw[purple] (8,3) circle (1pt);
   \filldraw[purple] (8,4) circle (1pt);
   \filldraw[purple] (9,0) circle (1pt);
   \filldraw[purple] (9,1) circle (1pt);
   \filldraw[purple] (9,3) circle (1pt);
   \filldraw[purple] (9,4) circle (1pt);
   \filldraw[purple] (9,5) circle (1pt);
   \filldraw[purple] (10,0) circle (1pt);
   \filldraw[purple] (10,1) circle (1pt);
   \filldraw[purple] (10,2) circle (1pt);
   \filldraw[purple] (10,3) circle (1pt);
   \filldraw[purple] (10,4) circle (1pt);
   \filldraw[purple] (10,5) circle (1pt);
   \filldraw[purple] (11,1) circle (1pt);
   \filldraw[purple] (11,2) circle (1pt);
   \filldraw[purple] (11,3) circle (1pt);
   \filldraw[purple] (11,4) circle (1pt);
   \filldraw[purple] (11,5) circle (1pt);
   \node at (16,4) {\underline{balancing condition}};
   \node at (16,3) {$1v_1+1v_2+1v_3=0$};
   \node at (16,2) {$1u_1+1u_2+2u_3=0$};
   \node at (16,1) {$2w_1+3w_2+1w_3=0$};
  \end{tikzpicture}
 \endpgfgraphicnamed
 \]
 \caption{A tropical curve in $\R^2$ and the balancing condition}
 \label{figure: tropical curve}
\end{figure}

The ties to classical algebraic geometry are given in terms of tropicalizations. Let $k$ be an algebraically closed field together with a valuation $v:k^\times\to\R$ with dense image. Let $X$ be a closed subvariety of the algebraic torus $(k^\times)^n$. Then the coordinatewise application of $v$ defines a subset of $\R^n$ whose topological closure is defined as the tropicalization $X^\trop$ of $X$. The structure theorem of tropical geometry asserts that $X^\trop$ can be subdivided into finitely many polyhedra and that the classical variety $X$ determines weights for the top dimensional polyhedra that satisfy the balancing condition. This can be extended to subvarieties of any toric variety, e.g.\ subvarieties of a projective space. 

\begin{figure}[ht]
 \[
 \resizebox{.9\linewidth}{!}{
 \includegraphics{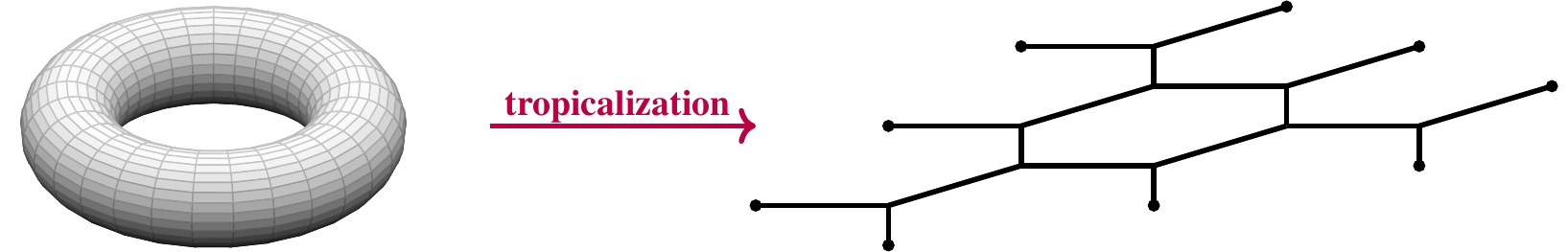}
 }
 \]
 \label{figure: tropicalization of an elliptic curve}
 \caption{Tropicalization of an elliptic curve, including its points at infinity}
\end{figure}

\medskip\noindent\textbf{The central problem.}
The definition of a tropical variety is problematic in several senses. First of all, the structure of a polyhedral complex of $X^\trop$ is not determined by the classical variety, but involves choices. Thus strictly speaking, the tropicalization of a classical variety is \emph{not} a tropical variety. Secondly, there is a discrepancy between the polynomial algebra over the tropical numbers and the functions determined by these polynomials: different polynomials can define the same function. These defects call for more sophisticated foundations of tropical geometry.

\subsubsection{Tropical scheme theory}

The above mentioned problem has been overcome by the paradigm changing approach of Jeff and Noah Giansiracusa (\cite{Giansiracusa-Giansiracusa16}), in which they realize the tropicalization $X^\trop$ of a classical variety $X$ as the set of $\T$-rational points of a $\T$-scheme $\cX$, which is a semiring scheme with a structure morphism to $\Spec\T$. This development can be seen as the tropical analogue of the passage from classical algebraic geometry, where varieties were considered as point sets in an ambient affine or projective space, to modern algebraic geometry in the sense of Grothendieck.

Soon after the initial paper \cite{Giansiracusa-Giansiracusa16} appeared, Maclagan and Rinc\'on (\cite{Maclagan-Rincon14}) have shown that the balancing condition can be recovered from the structure as a $\T$-scheme together with its embedding into an algebraic torus. 

To complete the step from embedded varieties to abstract schemes, the author has developed in \cite{L15b} a coordinate-free theory of tropicalizations, which is based on blueprints and blue schemes. In more detail, the Giansiracusa tropicalization of a classical variety is not only a semiring scheme, but comes with the richer structure as a blue scheme. The structure of the tropicalization as a blue scheme is sufficient to determine the weights of the underlying tropical variety. This makes it possible to pass from embedded tropical schemes to abstract tropical schemes. The gain of this change of perspective is that it applies, under suitable conditions, to more general situations, such as skeleta of Berkovich spaces and tropicalizations of moduli spaces of curves.

\subsubsection{Future applications}

At the time of writing, there are high hopes that this new approach to tropical geometry will lead to a realization of a conjectured tropical sheaf cohomology and subsequently allows for progress in Brill-Noether and Baker-Norine theory. It also might put tropical intersection theory on a new footing. Finally, we expect that tropical scheme theory will interplay with Connes and Consani's program around the Riemann hypothesis, which is based on idempotent semirings as tropical scheme theory is.

\begin{small}
 \bibliographystyle{alpha}
\newcommand{\etalchar}[1]{$^{#1}$}

\end{small}

\end{document}